\title{Representations of the Temperley-Lieb Algebra via a New Inner Product on Half-Diagrams}
\author{Emily Redelmeier}
\date{December, 2006}

\documentclass[11pt]{article}

\usepackage{graphicx}
\usepackage{graphics}
\usepackage{amsfonts}
\usepackage{amssymb}
\usepackage{amsthm}
\usepackage{amsmath}

\newcommand{\tr}{\mathop{\mathrm{tr}}\nolimits}

\newtheorem{thm}{Theorem}[section]
\newtheorem{lem}[thm]{Lemma}
\newtheorem{prop}[thm]{Proposition}
\newtheorem{cor}[thm]{Corollary}

\begin{document}

\maketitle

\section{The Temperley-Lieb Algebra}

The Temperley-Lieb algebra is an algebra given by generators and relations which arose in a problem in statistical mechanics considered by Temperley and Lieb \cite{MR0498284}.  Since then the algebra has been used in a wide variety of applications: the Witten-Reshetikhin-Turaev invariants of \(3\)-manifolds, Jones' analysis of the index of von Neumann subfactors \cite{MR696688}, Jones' polynomial invariant for knots, the four-colour problem (see, for example, \cite{MR1999738}), and Di Francesco's work on counting meanders in relation to Hilbert's sixteenth problem \cite{MR1608551,MR1462755}.

In \cite{MR696688} Jones calculated the irreducible representations of \(TL_{n}\left(q\right)\) for \(q\geq 2\) and for the semisimple part of \(TL_{n}\left(q\right)\) when \(q=2\cos\left(\frac{\pi}{n}\right)\).  This analysis was then used to show that the index of a subfactor was restricted to the union of a sequence and an interval.  This soon led to the polynomial invariant for knots found by Jones and now known as the Jones polynomial.

We shall work with a diagrammatic interpretation of the Temperley-Lieb algebra found by Kauffman \cite{MR899057} involving planar (or noncrossing) diagrams.

As these algebras are fundamental to many areas of mathematics and theoretical physics, detailed understanding of their representations is very important.  In this project we consider a space of half-diagrams and show that a natural inner product and an action of the Temperley-Lieb algebra may be given which for suitable values of a parameter realizes all of the irreducible representations.

The Temperley-Lieb algebra \(TL_{n}\left(q\right)\) or \(TL_{n}\), for some fixed integer \(n\) and some complex number \(q\), is defined as the algebra generated by \(1,e_{1},\ldots,e_{n-1}\) satisfying the {\em Jones relations}: \(e_{i}^{2}=qe_{i}\), \(e_{i}e_{j}=e_{j}e_{i}\) for \(\left|i-j\right|>1\), and \(e_{i}e_{i\pm 1}e_{i}=e_{i}\).  In this paper, \(q\) will always be a positive real number.  The Temperley-Lieb algebra can also be described by certain graphs and operations on those graphs, which we describe in this section.

A {\em pairing} of the integers \(\left\{1,\ldots,2n\right\}\) is a partition of these integers into sets of two.  We can interpret this as a graph on the integers \(\left\{1\,\ldots,2n\right\}\) with an edge connecting the two integers in the same set.

A {\em crossing} is a set of four integers \(i,j,k,l\in\left\{1,\ldots,2n\right\}\) such that \(i<k<j<l\), where \(i\) is paired with \(j\) and \(k\) with \(l\).  If a pairing has no crossings, it is said to be {\em noncrossing}.

This definition of crossing can be understood if paired integers (on the real line of the complex plane) are connected with arcs drawn in the upper half-plane.  It is also sometimes convenient to place the integers \(1\) through \(n\) from top to bottom along the left side of a rectangle and integers \(n+1\) through \(2n\) from bottom to top along the right, and draw curves connecting pairs inside the rectangle.  The same definition of a crossing still makes sense in this configuration.

\begin{figure}
\centering
\begin{tabular}{ccccc}
\includegraphics{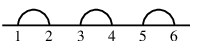}&\includegraphics{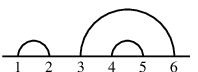}&\includegraphics{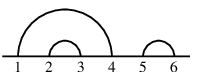}&\includegraphics{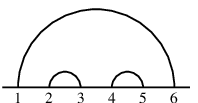}&\includegraphics{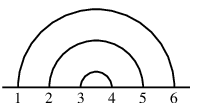}\\
\includegraphics{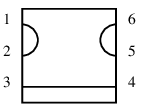}&\includegraphics{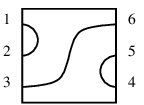}&\includegraphics{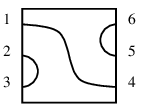}&\includegraphics{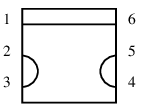}&\includegraphics{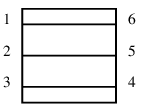}
\end{tabular}
\caption{The five noncrossing pairings on \(\left\{1,\ldots,6\right\}\), shown as arch diagrams and as strand diagrams}
\label{full-basis}
\end{figure}

We will often drop the adjective ``noncrossing'' and refer to noncrossing pairings simply as ``pairings''.  We will also refer to noncrossing pairings as ``diagrams'', or as ``arch diagrams'' (especially when we think of the points arranged on the real line) or ``strand diagrams'' (especially when we think of the points arranged in a rectangle).  There are five noncrossing pairings on \(\left\{1,\ldots,6\right\}\), shown in Figure \ref{full-basis}.

We can also define the Temperley-Lieb algebra \(TL_{n}\) as the vector space of formal linear combinations over the complex numbers of the noncrossing pairings on the integers \(\left\{1,\ldots,2n\right\}\).

\begin{figure}
\centering
\[\includegraphics{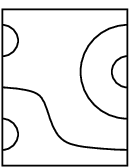}\times\includegraphics{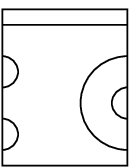}=\includegraphics{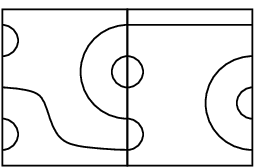}=q\includegraphics{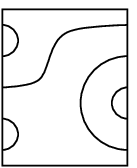}\]
\caption{Multiplication, represented diagrammatically}
\label{multiplication}
\end{figure}

We define multiplication as follows.  Any two basis elements \(e_{\mathfrak{a}}\) and \(e_{\mathfrak{b}}\) correspond to noncrossing strand diagrams \(\mathfrak{a}\) and \(\mathfrak{b}\).  We concatenate the strand diagrams by identifying the \(\left(2n-i+1\right)\)th point of the diagram \(\mathfrak{a}\) with the \(i\)th point of the diagram \(\mathfrak{b}\) for \(i\) with \(1\leq i\leq n\), as in Figure \ref{multiplication}.  In this concatenation, we get some number \(c\) of closed loops.  We can also construct a noncrossing diagram \(\mathfrak{c}\) from the concatenation by deleting all closed loops as well as the points that were identified in the concatenation.  (Whenever we speak of deleting a point, it has exactly two edges.  We identify these edges by replacing them with an edge connecting the two points originally connected to this point.)  We let the \(i\)th point of \(\mathfrak{c}\) be the \(i\)th point of \(\mathfrak{a}\) for \(1\leq i\leq n\) and the \(i\)th point of \(\mathfrak{b}\) for \(n+1\leq i\leq 2n\).  This corresponds to a pairing, since each remaining point has exactly one edge (as it had exactly one edge as a point of either \(\mathfrak{a}\) or \(\mathfrak{b}\) and was not identified with any other point), so it is connected to exactly one other point.  We then define the product \(e_{\mathfrak{a}}e_{\mathfrak{b}}=q^{c}e_{\mathfrak{c}}\), and extend this multiplication bilinearly over the entire vector space.  We can see from the diagrams that any new diagram produced this way is noncrossing, and the multiplication is associative.

\begin{figure}
\centering
\begin{tabular}{ccccc}
\includegraphics{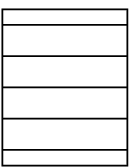}&\includegraphics{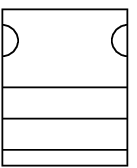}&\includegraphics{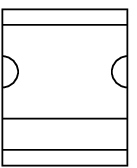}&\includegraphics{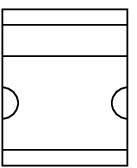}&\includegraphics{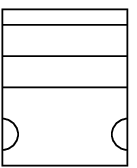}\\\(1\)&\(e_{1}\)&\(e_{2}\)&\(e_{3}\)&\(e_{4}\)
\end{tabular}
\caption{The generators of \(TL_{5}\)}
\label{generators}
\end{figure}

If we let \(\mathfrak{a}\) be the diagram in which the \(i\)th point is connected to the \(\left(2n-i+1\right)\)th point for all \(i\), \(1\leq i\leq n\), then \(e_{\mathfrak{a}}\), the basis element of \(TL_{n}\) associated with \(\mathfrak{a}\), is a right and left identity in this algebra, as we can see from Figure \ref{generators}.

We define \(e_{i}\), for each \(i\), \(1\leq i\leq n-1\), as the basis element of \(TL_{n}\) corresponding to the diagram in which the \(i\)th point is connected to the \(\left(i+1\right)\)th, the \(\left(2n-i+1\right)\)th point is connected to the \(\left(2n-i\right)\)th, and the \(j\)th point is connected to the \(\left(2n-j+1\right)\)th for all \(j\) not equal to \(i\) or \(i+1\).  We show \(e_{1},e_{2},e_{3},e_{4}\in TL_{5}\) in Figure \ref{generators}.

\begin{figure}
\centering
\[\begin{array}{c}\includegraphics{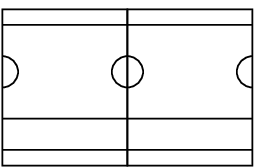}=q\includegraphics{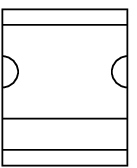}\\\includegraphics{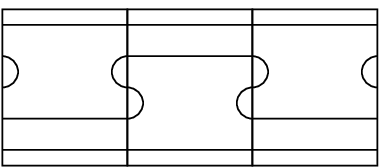}=\includegraphics{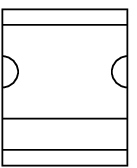}\\\includegraphics{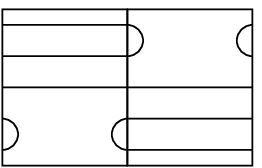}=\includegraphics{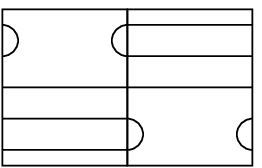}\end{array}\]
\caption{The Jones relations on diagrams}
\label{Jones}
\end{figure}

Thus the algebra we have defined using these graphs satisfies the Jones relations, which we can see in Figure \ref{Jones}.  In fact, these relations are sufficient to describe this algebra.  For a proof, see, for example, \cite{MR1343661}.

For any pairing \(\mathfrak{a}\), we define its transpose \(\mathfrak{a}^{t}\) as the pairing we get by relabelling the \(i\)th point \(2n-i+1\), for \(1\leq i\leq 2n\).  (Since \(i\mapsto 2n-i+1\) is a bijection on \(\left\{1,\ldots,2n\right\}\), this \(\mathfrak{a}^{t}\) is a pairing, and since it reverses the ordering, it does not produce any crossings.)  As \(\left(\mathfrak{a}^{t}\right)^{t}=\mathfrak{a}\), we can define an involution \(*\) on \(TL_{n}\) by letting \(e_{\mathfrak{a}}^{*}=e_{\mathfrak{a}^{t}}\) (where \(e_{\mathfrak{a}}\) is the basis element associated with diagram \(\mathfrak{a}\)) and extending conjugate linearly to the whole space.  If we renumber the points of the concatenation of diagrams \(\mathfrak{a}\) and \(\mathfrak{b}\) in this way, we can see that this gives us the concatenation of \(\mathfrak{b}^{t}\) and \(\mathfrak{a}^{t}\), so \(\left(e_{\mathfrak{a}}e_{\mathfrak{b}}\right)^{*}=e_{\mathfrak{b}}^{*}e_{\mathfrak{a}}^{*}\), so the involution reverses the order of multiplication.

We define the {\em closure} of a pairing \(\mathfrak{a}\) by identifying the \(i\)th point with the \(\left(2n-i+1\right)\)th point.  Since each point in \(\mathfrak{a}\) has one edge, and is identified with exactly one other point, each point in the closure has two edges, and hence this diagram consists of some number \(c\) of closed loops.  We define the trace of the associated basis element \(e_{\mathfrak{a}}\) as \(\tr\left(e_{\mathfrak{a}}\right)=q^{c}\), and extend this definition linearly to the rest of \(TL_{n}\).

As we now have a linear trace function, a bilinear multiplication and a conjugate linear involution on \(TL_{n}\), we can define a sesquilinear function \(\langle e,f\rangle=\tr\left(ef^{*}\right)\).  We show later that when the first \(n\) Chebyshev polynomials of the second kind are positive at \(q\), this function is positive and hence an inner product, but that it does not generally satisfy the positivity requirement for an inner product for other \(q\).  However, we will refer to it as an inner product.

\section{Half-Diagrams}

If we consider the left half of a pairing, that is, the points \(\left\{1,\ldots,n\right\}\) of a pairing on \(2n\) points, \(2p\) of the points are paired with another point in \(\left\{1,\ldots,n\right\}\), for some \(p\).  The other \(n-2p\) of the points, which were connected to a point on the opposite side (that is, one of the points of \(\left\{n+1,\ldots,2n\right\}\)), can now be thought of as being connected to a point at infinity, or as having a {\em through-string}.  We will also refer to a pair of points in a full diagram, one of which is in \(\left\{1,\ldots,n\right\}\) and one of which is in \(\left\{n+1,\ldots,2n\right\}\), as a through-string.

We define a {\em half-diagram} as a partition of the points \(\left\{1,\ldots,n\right\}\) into sets of one or two each, where the sets of two are pairings and the sets of one are thought to have a through-string.  We extend the definition of a crossing to include sets of three integers \(i\), \(j\) and \(k\) such that \(i<k<j\), where \(i\) is paired with \(j\) and \(k\) has a through-string.  Again, this definition can be understood if the paired integers are connected by arches in the upper half-plane and the integers with through-strings are connected upward to a point at infinity.

\begin{figure}
\centering
\includegraphics{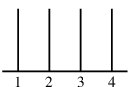}\\\includegraphics{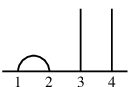}\includegraphics{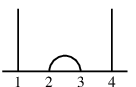}\includegraphics{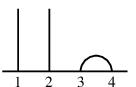}\\\includegraphics{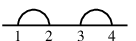}\includegraphics{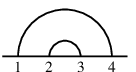}
\caption{The six noncrossing half-diagrams on four points}
\label{half-basis}
\end{figure}

Again, we will generally drop the adjective ``noncrossing'' and refer to noncrossing half diagrams as half diagrams or half pairings.  There are six noncrossing half-diagrams on four points: one with no pairs and four through-strings, three with one pair and two through-strings, and two with two pairs and no through-strings, all shown in Figure \ref{half-basis}.

We can again construct a vector space of formal linear combinations of the noncrossing half-diagrams on \(n\) points.  We denote this \(U\left(n\right)\), and we denote the subspace generated by the diagrams with \(p\) pairs and \(n-2p\) through-strings by \(U\left(n;p\right)\).  We will use Latin letters to denote elements of \(TL_{n}\) and Greek letters to denote elements of \(U\left(n\right)\).

\begin{figure}
\centering
\[\begin{array}{l}\includegraphics{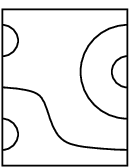}\cdot\includegraphics{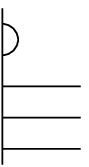}=\includegraphics{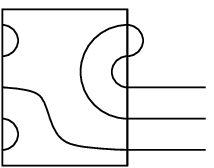}=0\\\includegraphics{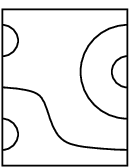}\cdot\includegraphics{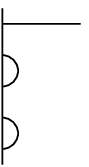}=\includegraphics{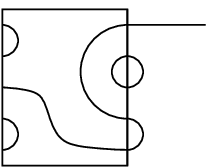}=q\includegraphics{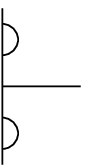}\end{array}\]
\caption{The action of a basis element of \(TL_{n}\) on two elements of \(U\left(n\right)\)}
\label{action}
\end{figure}

We define the action of \(TL_{n}\) on \(U\left(n\right)\) by defining the action of the basis elements of \(TL_{n}\) on the basis elements of \(U\left(n\right)\).  If \(\mathfrak{p}\) is a noncrossing pairing on \(2n\) points and \(\mathfrak{a}\) is a noncrossing half-pairing on \(n\) points, we concatenate their diagrams by identifying the \(i\)th point of \(\mathfrak{a}\) with the \(\left(2n-i+1\right)\)th point of \(\mathfrak{p}\).  If there is an edge path between two through-strings of \(\mathfrak{a}\), we define \(e_{\mathfrak{p}}\xi_{\mathfrak{a}}=0\).  Otherwise, this produces some number \(c\) of closed loops.  We can define a new noncrossing half diagram \(\mathfrak{b}\) by deleting all closed loops from the concatenation of \(\mathfrak{p}\) and \(\mathfrak{a}\), as well as all of the identified points, and letting the points \(\left\{1,\ldots,n\right\}\) of \(\mathfrak{b}\) be points \(\left\{1,\ldots,n\right\}\) of \(\mathfrak{p}\).  As each of the points aside from those in \(\left\{1,\ldots,n\right\}\) has been identified with one other point, each of these points has two edges.  So each of the points in \(\left\{1,\ldots,n\right\}\) is eventually connected by a sequence of edges to either another point in \(\left\{1,\ldots,n\right\}\) or to the point at infinity.  So if we delete the identified points, we are left with a half-diagram \(\mathfrak{b}\) which we can see is noncrossing.  We define \(e_{\mathfrak{p}}\xi_{\mathfrak{a}}=q^{c}\xi_{\mathfrak{b}}\).  We show some examples in Figure \ref{action}.  We can then extend this action linearly in each of \(TL_{n}\) and \(U\left(n\right)\).  As the concatenation process respects multiplication in \(TL_{n}\), this gives us a representation of \(TL_{n}\) on \(U\left(n\right)\).

If \(e_{\mathfrak{p}}\xi_{\mathfrak{a}}\neq 0\), then no through-string is connected to another through-string in the concatenation of \(\mathfrak{p}\) with \(\mathfrak{a}\).  Then if \(\mathfrak{a}\) has \(n-2p\) through-strings, each eventually connects to one of the points in \(\left\{1,\ldots,n\right\}\).  So \(\mathfrak{b}\) has \(n-2p\) through-strings as well.  Thus this representation can be reduced to representations on \(U\left(n;p\right)\).

Since the noncrossing half diagrams on \(2n\) points with no through-strings are the noncrossing diagrams on \(2n\) points, we can think of \(U\left(2n;n\right)\) as \(TL_{n}\).

We can also think of a noncrossing diagram \(\mathfrak{p}\) on \(2n\) points as a pair of noncrossing half-diagrams on \(n\) points.  The pairings on the points \(1\) through \(n\) define one half-diagram, where the points connected to points outside of this set now have a through-string.  Since the other points \(l\) are all greater than \(n\), there can be no \(k\) with \(i<k<j\) where \(k\) now has a through-string, so this new diagram is noncrossing.  Likewise, we can construct a half-diagram from the points \(\left\{n+1,\ldots,2n\right\}\) of \(\mathfrak{p}\) by relabelling the \(\left(2n-i+1\right)\)th point \(i\), where again a point connected to a point outside this set now has a through-string, and likewise this half-diagram will be noncrossing.

Conversely, given two noncrossing half-diagrams \(\mathfrak{a}\) and \(\mathfrak{b}\) with the same number \(n-2p\) of through-strings, we can construct a unique noncrossing diagram \(\mathfrak{p}\).  We let the \(i\)th point of \(\mathfrak{p}\) be the \(i\)th point of \(\mathfrak{a}\), and the \(\left(2n-i+1\right)\)th point be the \(i\)th point of \(\mathfrak{b}\).  \(\mathfrak{a}\) has through-strings at \(i_{1},\ldots,i_{2n-p}\) with \(i_{1}<\ldots<i_{2n-p}\), and \(\mathfrak{b}\) has through-strings at \(j_{1},\ldots,j_{2n-p}\) with \(j_{1}<\ldots<j_{2n-p}\).  If we connect \(i_{t}\) to \(2n-j_{t}+1\) for \(1\leq t\leq k\) in \(\mathfrak{p}\), we have a pairing on \(\left\{1,\ldots,2n\right\}\).  These new arches cannot cross any old arches, since there are no paired \(i\) and \(j\) such that \(i<i_{1}<j\) or \(i<j_{t}<j\).  Furthermore, for \(t_{1}<t_{2}\), \(i_{t_{1}}<i_{t_{2}}<2n-j_{t_{2}}+1<2n-j_{t_{1}}\), so the new arches do not cross each other.  Thus \(\mathfrak{p}\) is a noncrossing pairing.

Furthermore, if we pair the \(i_{t}\) and the \(j_{t}\) according to any other permutation (that is, \(i_{t}\) is paired with \(j_{\sigma\left(t\right)}\) for some \(\sigma\) not the identity), we must get a crossing, since there must be some \(t_{1}\) such that \(t_{1}<\sigma\left(t_{1}\right)\) (since there must be at least one \(t\) such that \(t\neq\sigma\left(t\right)\), and if there is one \(t\) such that \(t>\sigma\left(t\right)\), then at least one of the integers less than or equal to \(\sigma\left(t\right)\) must be mapped to at integer greater than \(\sigma\left(t\right)\)).  Then there must be a \(t_{2}\) such that \(t_{1}<t_{2}\) but \(\sigma\left(t_{2}\right)<\sigma\left(t_{1}\right)\) (since at least one of the integers greater than \(t_{1}\) must be mapped to an integer less than \(\sigma\left(t_{1}\right)\)).  Then \(i_{t_{1}}<i_{t_{2}}<2n-j_{\sigma\left(t_{1}\right)}+1<2n-j_{\sigma\left(t_{2}\right)}+1\), a crossing.  So the noncrossing pairing \(\mathfrak{p}\) associated to the two half-pairings \(\mathfrak{a}\) and \(\mathfrak{b}\) is unique.

Thus, \(TL_{n}\cong\left(U\left(n,0\right)\otimes\left(n,0\right)\right)\oplus\ldots\oplus\left(U\left(n,\lfloor\frac{n}{2}\rfloor\right)\right)\), where \(\lfloor\frac{n}{2}\rfloor\) is the greatest integer less than or equal to \(\frac{n}{2}\), and \(\xi_{\mathfrak{a}}\otimes\xi_{\mathfrak{b}}=e_{\mathfrak{p}}\), for \(\mathfrak{a}\), \(\mathfrak{b}\) and \(\mathfrak{p}\) as above.  We can then abuse this notation write \(\mathfrak{p}=\mathfrak{a}\otimes\mathfrak{b}\).

\begin{figure}
\centering
\[\begin{array}{c}\langle\includegraphics{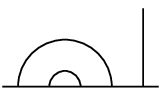},\includegraphics{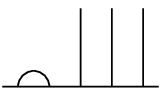}\rangle=\includegraphics{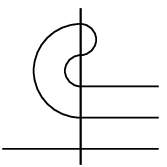}=0\\\langle\includegraphics{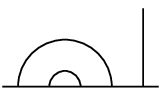},\includegraphics{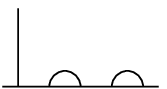}\rangle=\includegraphics{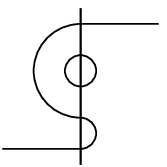}=q\end{array}\]
\caption{Two examples of inner products}
\label{inner}
\end{figure}

We define a sesquilinear function on \(U\left(n\right)\).  For \(\mathfrak{a},\mathfrak{b}\in U\left(n\right)\), we identify the \(i\)th point of \(\mathfrak{a}\) with the \(i\)th point of \(\mathfrak{b}\).  If we get any edge path that connects a through-string to another through-string on the same side, we let \(\langle\xi_{\mathfrak{a}},\xi_{\mathfrak{b}}\rangle=0\).  Otherwise, we get some number \(c\) of closed loops, and we let \(\langle\xi_{\mathfrak{a}},\xi_{\mathfrak{b}}\rangle=q^{c}\).  We can think of edge paths connecting through-strings on the same edge as contributing a factor of \(0\), edge paths connect a through-string on one side with a through-string on the other as contributing a factor of \(1\), and closed loops as contributing a factor of \(q\).  We extend this sesquilinearly to all of \(U\left(n\right)\).  Two examples are shown in Figure \ref{inner}.

This definition of inner product agrees on \(U\left(2n;n\right)\) with the previously defined inner product on \(TL_{n}\).  When we multiply \(\mathfrak{p}\) by \(\mathfrak{q}^{t}\), we identify the \(\left(2n-i+1\right)\)th point of \(\mathfrak{p}\) with the \(i\)th point of \(\mathfrak{q}^{t}\), that is, the \(\left(2n-i+1\right)\)th point of \(\mathfrak{q}\), for \(1\leq i\leq n\).  When we take the closure of \(\mathfrak{p}\mathfrak{q}^{t}\), we identify the \(i\)th point of the diagram with the \(\left(2n-i+1\right)\)th, that is, the \(i\)th point of \(\mathfrak{p}\) with the \(\left(2n-i+1\right)\)th point of \(\mathfrak{q}^{t}\), that is, the \(i\)th point of \(\mathfrak{q}\), for \(1\leq i\leq n\).  Thus, we have identified the \(i\)th point of \(\mathfrak{p}\) with the \(i\)th point of \(\mathfrak{q}\) for all \(i\), \(1\leq i\leq 2n\), as we do when we take the inner product of half-diagrams.  As there are no through-strings, in both cases, the inner product is \(q^{c}\), where \(c\) is the number of closed loops in the resulting diagram, or the trace of \(\mathfrak{p}\mathfrak{q}^{t}\).

Again, if the first \(n\) Chebyshev polynomials of the second kind are positive at \(q\), this function is also positive, and is hence an inner product.  Again, this is not the case for general \(q\), but we will refer to it as an inner product nonetheless.

Furthermore, the involution \(*\) on \(TL_{n}\), which flips a diagram left-to-right, is an adjoint.  In fact, we will show below that this is the unique sesquilinear function with this property, up to scalar multiples.

\section{Indexing the Bases}

We will show that the number of noncrossing diagrams on \(2n\) points is \(c_{n}=\frac{1}{n+1}\binom{2n}{n}\), the \(n\)th {\em Catalan number}.  The number of noncrossing half diagrams on \(n\) points with \(p\) pairs of points is \(c_{n,n-2p}=\binom{n}{p}-\binom{n}{p-1}=\binom{n-1}{p}-\binom{n-1}{p-1}\), called a {\em generalized Catalan number} by Di Franceso \cite{MR1608551}.  Quite a few mathematical objects are counted by the Catalan numbers, and many of these objects can be generalized to an object counted by the generalized Catalan numbers.  Several of these are useful in this paper.

\subsection{Catalan Constructions}

\begin{figure}
\centering
\includegraphics{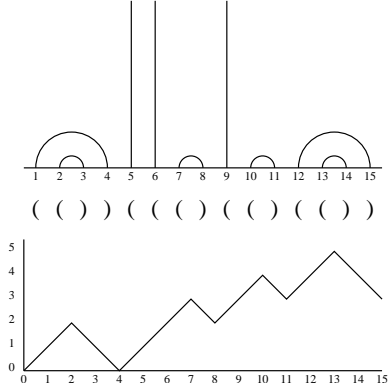}
\caption{The equivalent half-pairing, generalized bracket diagram, and generalized Dyck path, corresponding to generalized restricted sequence \(\left(2,1,3,4,5,4\right)\)}
\label{Catalan}
\end{figure}

Noncrossing arch diagrams, bracket diagrams, Dyck paths and restricted sequences, all of which will be defined below, are all counted by the Catalan numbers, and each has a generalization counted by the generalized Catalan numbers.  We will generally call these objects the generalized version (that is, generalized bracket diagram, generalized Dyck path, etc.) when the distinction has to be made between the generalized object and the original, although we will usually refer to them as bracket diagrams, Dyck paths, etc. as well.  In Figure \ref{Catalan} we show an example of a half-diagram along with its equivalent bracket diagram, Dyck path and restricted sequence.

There are obvious bijections between these sets, which we will describe below.  We will generally use the same symbol to refer to the arch diagram as well as the corresponding bracket diagram, Dyck path and restricted sequence.  While the basis element of \(TL_{n}\) corresponding to an arch diagram \(\mathfrak{p}\) can be thought of as being that arch diagram, we will always denote it \(e_{\mathfrak{p}}\).  This way a second basis, whose elements we will denote \(e_{\mathfrak{p}}^{\prime}\), can also be indexed by the same set.  Likewise, we will use the same symbol, say, \(\mathfrak{a}\), to refer to a half diagram along with the corresponding generalized bracket diagram, generalized Dyck path and generalized restricted sequence, and use \(\xi_{\mathfrak{a}}\) to denote the corresponding basis element of \(U\left(n\right)\) and \(\xi_{\mathfrak{a}}^{\prime}\) to denote the element of the second basis.

\subsubsection{Bracket Diagrams}

A {\em bracket diagram} is a sequence of \(n\) opening brackets and \(n\) closing brackets such that, among the first \(i\) brackets, the number of closing brackets does not exceed the number of closing brackets, for all \(i\), \(1\leq i\leq 2n\).  This corresponds to our notion of a legitimate bracket diagram.  A {\em generalized bracket diagram} is a sequence of \(n\) brackets such that, again, the number of closing brackets among the first \(i\) brackets does not exceed the number of opening brackets, although we do not require there to be the same number of opening brackets as closing brackets overall.  Westbury's  definition of a bracket diagram in \cite{MR1343661} is equivalent to our definition of a generalized bracket diagram; however, Westbury's diagrams have dots instead of our unpaired opening brackets.  As we will show below, which opening brackets are paired and which are unpaired is unambiguous, so these constructions are equivalent.

We define a height function \(h_{i}\left(\mathfrak{a}\right)\), equal to the number of opening brackets minus the number of closing brackets among the first \(i\) brackets of \(\mathfrak{a}\).  This function takes on nonnegative integers, and may only change by \(\pm 1\) when \(i\) is increased or decreased by \(1\).

Given a noncrossing half-diagram \(\mathfrak{a}\), we can construct a generalized bracket diagram.  If the \(i\)th point is connected to the \(j\)th point where \(i<j\), or there is a through-string at the \(i\)th point, we let the \(i\)th bracket be an opening bracket.  Otherwise (that is, if the \(i\)th point is connected to the \(j\)th point, where \(j<i\)), we let the \(i\)th bracket be a closing bracket.  For any \(i\), among the first \(i\) points of \(\mathfrak{a}\), each point connected to a point of smaller index must be connected to another point among these first \(i\) points, so among the first \(i\) brackets, the number of closing brackets will never exceed the number of opening brackets.

Conversely, given a bracket diagram, we can construct a noncrossing half-pairing by pairing the points according to the pairings of the bracket.  In this way, each opening bracket is either unpaired or corresponds to a point paired with a point of greater index, and each closing bracket is paired with a bracket of smaller index.  Furthermore, we have no crossings.

More formally, we pair each \(j\), where the \(j\)th bracket is a closing bracket, with the largest \(i\) such that \(i<j\) and \(h_{i}\left(\mathfrak{a}\right)=h_{j-1}\left(\mathfrak{a}\right)\).

We know that such an \(i\) must exist.  The \(j\)th bracket is a closing bracket, and hence \(j>1\) (or the number of closing brackets in the first \(j\) brackets would exceed the number of opening brackets), so there is a \(\left(j-1\right)\)th bracket.  \(h_{j-1}\left(\mathfrak{a}\right)-h_{i}\left(\mathfrak{a}\right)=0\) for \(i=j-1\), although the \(i\)th bracket might not be an opening bracket.  However, if there is an \(i<j\) where \(h_{j-1}\left(\mathfrak{a}\right)-h_{i}\left(\mathfrak{a}\right)=0\) and the \(i\)th bracket is a closing bracket, then there is one fewer closing bracket in the interval \(\left[1,i-1\right]\) than in \(\left[1,i\right]\).  Then \(h_{j-1}\left(\mathfrak{a}\right)-h_{i-1}\left(\mathfrak{a}\right)=h_{j-1}\left(\mathfrak{a}\right)-h_{i}\left(\mathfrak{a}\right)-1=-1\).  \(h_{i}\left(\mathfrak{a}\right)-h_{j-1}\left(\mathfrak{a}\right)\), as a function of \(i\), can only change by \(\pm 1\) when \(i\) is decreased by \(1\), and \(h_{j-1}\left(\mathfrak{a}\right)-h_{0}\left(\mathfrak{a}\right)=h_{j-1}\left(\mathfrak{a}\right)\geq 0\), so there must be an \(i^{\prime}\) with \(i^{\prime}<i\) such that \(h_{j-1}\left(\mathfrak{a}\right)-h_{i^{\prime}}\left(\mathfrak{a}\right)=0\).  As there are only finitely many possible indices, there must be an \(i\) with \(h_{j-1}\left(\mathfrak{a}\right)-h_{i}\left(\mathfrak{a}\right)\) which does not satisfy the hypothesis that the \(i\)th bracket is a closing bracket.  Thus this process is well defined.

We can also see that \(h_{j-1}\left(\mathfrak{a}\right)-h_{k}\left(\mathfrak{a}\right)\geq 0\) for any \(k\), \(i<k<j\), since otherwise \(h_{j-1}\left(\mathfrak{a}\right)-h_{l}\left(\mathfrak{a}\right)\) (as a function of \(l\)) must decrease from \(0\) to \(-1\) as \(l\) is decreased by \(1\).  Then there must be an \(l>i\) such that \(h_{j-1}\left(\mathfrak{a}\right)-h_{l}\left(\mathfrak{a}\right)=0\) and such that the \(l\)th bracket is an opening bracket.

If \(j\) is paired with \(i\) by this rule, and \(l>j\) where the \(l\)th bracket is a closing bracket, and \(h_{l-1}\left(\mathfrak{a}\right)\geq h_{i}\left(\mathfrak{a}\right)\), then by an argument like that above, \(h_{l-1}\left(\mathfrak{a}\right)-h_{j-1}\left(\mathfrak{a}\right)\geq 0\) so \(h_{l-1}\left(\mathfrak{a}\right)-h_{j}\left(\mathfrak{a}\right)>0\), but \(h_{l-1}\left(\mathfrak{a}\right)-h_{l-1}\left(\mathfrak{a}\right)=0\).  Given a \(k\), \(j<k<l\) such that \(h_{l-1}\left(\mathfrak{a}\right)-h_{k}\left(\mathfrak{a}\right)=0\) and the \(k\)th bracket is a closing bracket, \(h_{l-1}\left(\mathfrak{a}\right)-h_{k-1}\left(\mathfrak{a}\right)=-1\), so we can find a \(k^{\prime}\) with \(j<k^{\prime}<k\) such that \(h_{l-1}\left(\mathfrak{a}\right)-h_{k^{\prime}}\left(\mathfrak{a}\right)=0\).  Again, as there are only finitely many possible indices, there must be a \(k\), \(j<k<l\), such that \(h_{l-1}\left(\mathfrak{a}\right)=h_{k}\left(\mathfrak{a}\right)=0\) and the \(k\)th bracket is an opening bracket.  As \(k>i\), \(l\) is paired with a point whose index is greater than \(j\), and hence greater than \(i\), under this rule.  Thus each closing bracket is paired with a different opening bracket, and hence this is legitimately a pairing.

Furthermore, if \(j\) is paired with \(i\) under this rule, and \(i<k<j\), then \(h_{k}\left(\mathfrak{a}\right)\geq h_{i}\left(\mathfrak{a}\right)\).  If there is an \(l\) with \(l>j\) such that the \(l\)th bracket is a closing bracket and \(h_{l-1}\left(\mathfrak{a}\right)=h_{k}\left(\mathfrak{a}\right)\), then as shown above, the \(l\)th point is paired with the \(m\)th point for some \(m\) with \(j<m<l\).  So the \(k\)th point is either unpaired or paired with a point between \(k\) and \(j\).  Thus a crossing cannot occur.

We now show that this is the only noncrossing half-diagram which corresponds to the given bracket diagram.  If a bracket diagram \(\mathfrak{a}\) corresponds to a pairing, then each point \(j\) where the \(j\)th point is a closing bracket must be paired with an \(i\), \(i<j\).  The \(i\)th bracket must be an opening bracket, since it is paired with a point of greater index.  Since all the points in the interval \(\left(i,j\right)\) must be paired with other points in the interval (or we have a crossing), there must be the same number of opening brackets as closing brackets in \(\left(i,j\right)\), so \(h_{i}\left(\mathfrak{a}\right)=h_{j-1}\left(\mathfrak{a}\right)\).  Finally, if there is a \(k\), \(i<k<j\) such that \(h_{k}\left(\mathfrak{a}\right)=h_{j-1}\left(\mathfrak{a}\right)\) and the \(k\)th bracket is an opening bracket, then there must be a largest such \(k\).  The \(k\)th point cannot have a through-string, so it must be paired with an \(l\), where the \(l\)th bracket is a closing bracket.  Then \(h_{l-1}\left(\mathfrak{a}\right)=h_{k}\left(\mathfrak{a}\right)=h_{j-1}\left(\mathfrak{a}\right)\).  We know that we can construct a legitimate noncrossing pairing by the original rule, in which the \(j\)th point would be connected to the \(k\)th point, and the \(l\)th point must then be connected to the \(m\)th point for some \(m\), \(k<m<l\).  Then the \(m\)th bracket is an opening bracket and \(h_{m}\left(\mathfrak{a}\right)=h_{l-1}\left(\mathfrak{a}\right)=h_{j-1}\left(\mathfrak{a}\right)\).  This contradicts our assumption that \(k\) is the largest index for with the \(k\)th bracket is an opening bracket and \(h_{k}\left(\mathfrak{a}\right)=h_{j-1}\left(\mathfrak{a}\right)\).  So we cannot construct an alternate noncrossing half-pairing corresponding to this bracket diagram.

Thus we have a bijection between half-diagrams on \(n\) points and bracket diagrams with \(n\) brackets.  In particular, if a half-diagram has \(p\) arches, then the corresponding bracket diagram will have \(n-p\) opening brackets and \(p\) closing brackets, and \(h_{n}\left(\mathfrak{a}\right)=n-2p\).

\subsubsection{Dyck paths}

A {\em Dyck path} is a sequence of points in \(\mathbb{Z}^{2}\) from \(\left(0,0\right)\) to \(\left(0,2n\right)\) such that if the \(i\)th point is \(\left(x,y\right)\), the \(\left(i+1\right)\)th point is either \(\left(x+1,y+1\right)\) or \(\left(x+1,y-1\right)\) and such that each point has a nonnegative \(y\) coordinate.  A {\em generalized Dyck path} is a sequence of points in \(\mathbb{Z}^{2}\) from \(\left(0,0\right)\) to \(\left(n,n-2p\right)\) for some \(p\) such that, again, if the \(i\)th point is \(\left(x,y\right)\), the \(\left(i-1\right)\)th point is either \(\left(x+1,y+1\right)\) or \(\left(x+1,y-1\right)\) and the \(y\) coordinate is never negative.  We will think of \(\left(0,0\right)\) as the zeroeth point, so the \(i\)th point always has \(x\) coordinate \(i\).

Given a generalized bracket diagram \(\mathfrak{a}\), we can construct a generalized Dyck path by letting the \(i\)th point be \(\left(x+1,y+1\right)\) (where the \(\left(i-1\right)\)th point is \(\left(x,y\right)\)) when the \(i\)th bracket is an opening bracket and \(\left(x+1,y-1\right)\) when the \(i\)th bracket is a closing bracket.  We can see that \(h_{i}\left(\mathfrak{a}\right)\) is the \(y\) coordinate of the \(i\)th point.  Then the \(y\) coordinate of a point is never negative, so this is a generalized Dyck path.

Conversely, given a Dyck path \(\mathfrak{a}\), we can construct a bracket diagram by letting the \(i\)th bracket be an opening bracket when the \(i\)th point is \(\left(x+1,y+1\right)\) and a closing bracket when the \(i\)th point is \(\left(x+1,y-1\right)\).  As the \(y\) coordinate of the \(i\)th point is nonnegative, there must have been at least as many \(\left(x+1,y+1\right)\) steps as \(\left(x+1,y-1\right)\), and hence there are at least as many opening brackets as closing brackets among the first \(i\) brackets.  So this is a legitimate bracket diagram.

We thus have a bijection between generalized bracket diagrams and generalized Dyck paths, and hence between half-diagrams and generalized Dyck paths.  We can think of opening brackets as upward steps and closing brackets as downward steps.

We say that \(\mathfrak{a}\) has a {\em maximum} at \(i\) if the \(y\) coordinate of the \(i\)th point is greater than the \(y\) coordinate of both the \(\left(i-1\right)\)th point and the \(\left(i+1\right)\)th point (so the \(i\)th bracket of \(\mathfrak{a}\) is an opening bracket and the \(\left(i+1\right)\)th bracket is a closing bracket).  The path shown in Figure \ref{Catalan} has maxima at \(2\), \(7\), \(10\) and \(13\).  We say that \(\mathfrak{a}\) has a {\em minimum} at \(i\) if the \(y\) coordinate of the \(i\)th point is smaller than the \(y\) coordinate of the \(\left(i-1\right)\)th point and the \(\left(i+1\right)\)th point (so the \(i\)th bracket is a closing bracket and the \(\left(i+1\right)\)th bracket is an opening bracket).  The path in Figure \ref{Catalan} has minima at \(4\), \(8\) and \(11\).  We say that \(\mathfrak{a}\) has a {\em slope} at \(i\) if the \(y\) coordinate of the \(i\)th point is between the \(y\) coordinate of the \(\left(i-1\right)\)th point and the \(y\) coordinate of the \(\left(i+1\right)\)th point (so the \(i\)th bracket and the \(\left(i+1\right)\)th bracket are either both opening brackets or both closing brackets).  The path in Figure \ref{Catalan} has slopes at \(1\), \(3\), \(5\), \(6\), \(9\), \(12\) and \(14\).  If \(h_{i-1}\left(\mathfrak{a}\right)<h_{i}\left(\mathfrak{a}\right)<h_{i+1}\left(\mathfrak{a}\right)\), so the \(i\)th bracket and the \(\left(i+1\right)\)th bracket are opening brackets, then we say \(\mathfrak{a}\) has an {\em increasing slope} at \(i\), and if \(h_{i-1}\left(\mathfrak{a}\right)>h_{i}\left(\mathfrak{a}\right)>h_{i+1}\left(\mathfrak{a}\right)\), so the \(i\)th bracket and the \(\left(i+1\right)\)th bracket are both closing brackets, we say that \(\mathfrak{a}\) has a {\em decreasing slope} at \(i\).  The slopes at \(1\), \(5\), \(6\), \(9\) and \(12\) are increasing and those at \(3\) and \(14\) are decreasing.

If there are \(n-2p\) through-strings in a half-diagram \(\mathfrak{a}\), then there are \(n-2p\) more opening brackets than closing brackets in the corresponding bracket diagram, so \(h_{n}\left(\mathfrak{a}\right)=n-2p\), and hence the final point of the corresponding generalized Dyck path is \(\left(n,n-2p\right)\).  So the half-diagrams with \(n-2p\) through-strings correspond to the generalized Dyck paths from \(\left(0,0\right)\) to \(\left(n,n-2p\right)\), and the full diagrams correspond to the original Dyck paths.

\subsubsection{Restricted Sequences}

A {\em restricted sequence}, defined by Ko and Smolinsky in \cite{MR1105701} as a description of a noncrossing pairing, is a sequence of \(n\) positive integers \(\left(a_{1},\ldots,a_{n}\right)\) such that \(a_{n}=1\) (if the sequence is nonempty) and, for each \(i\) with \(1\leq i<n\), \(a_{i+1}\geq a_{i}-1\).  We define a generalized restricted sequence as a sequence of \(p\) positive integers \(\left(a_{1},\ldots,a_{p}\right)\) (where \(p\leq\lfloor\frac{n}{2}\rfloor\)) where \(a_{i+1}\geq a_{i}-1\) for all \(i\) with \(1\leq i<p\), but we relax the condition that the last integer be \(1\).  Usually the number of points \(n\) in the associated half-diagram will be clear from the context.

We associate the empty sequence with the diagram on \(n\) points in which each point has a through-string.  As this is the only diagram on \(0\) points or \(1\) point, the diagrams on \(0\) points or \(1\) point are uniquely described by the generalized restricted sequences with \(p\) integers, \(p\leq\frac{n}{2}=0\).  Assume inductively that half-diagrams on \(n\) points are described uniquely by generalized restricted sequences with \(p\) integers, for each \(p\leq\frac{n}{2}\).  A noncrossing pairing \(\mathfrak{a}\) on \(n+2\) points with at least one arch must have a unique leftmost innermost arch, that is, an arch connecting the \(a_{1}\)th point to the \(\left(a_{1}+1\right)\)th point, where we choose \(a_{1}\) to be as small as possible.  The first closing bracket in the corresponding bracket diagram must belong to this arch, since the first bracket cannot be a closing bracket, and hence the first closing bracket must follow an opening bracket, which together form such an arch.  If we remove this arch and its points (and consider the \(j\)th point to the \(\left(j-2\right)\)th point of the new diagram for every \(j>i\)), we have a half-diagram \(\mathfrak{b}\) on \(n\) points, which, according to our inductive hypothesis, corresponds to a unique restricted sequence \(\left(b_{1},\ldots,b_{p}\right)\) for some \(p\leq\frac{n-2}{2}\).  Then we let \(\mathfrak{a}\) correspond to the sequence \(\left(a_{1},b_{1},\ldots,b_{p}\right)\).  The leftmost innermost arch of \(\mathfrak{b}\) connects the \(b_{1}\)th point with the \(\left(b_{1}+1\right)\)th point.  If \(a_{1}>b_{1}+1\), then the points \(b_{1}\) and \(b_{1}+1\) in \(\mathfrak{b}\) correspond to points \(b_{1}\) and \(b_{1}+1\) in \(\mathfrak{a}\), and hence \(\mathfrak{a}\) has an arch connecting \(b_{1}\) and \(b_{1}+1\), and \(b_{1}<a_{1}\).  This contradicts our assumption that \(a_{1}\) was the smallest integer for which \(\mathfrak{a}\) has an arch joining the \(a_{1}\)th point to the \(\left(a_{1}+1\right)\)th point.  So \(b_{1}\geq a_{1}-1\).  As \(b_{i+1}\geq b_{i}-1\) for all \(i\) with \(1\leq i<p\), \(\left(a_{1},b_{1},\ldots,b_{p}\right)\) is a legitimate restricted sequence, and \(p+1\leq\frac{n+2}{2}\).  As \(a_{1}\) is unique, \(\mathfrak{a}\) corresponds to a unique generalized restricted sequence on \(p+1\) integers with \(p\leq\frac{n+2}{2}\).

Conversely, given a restricted sequence \(\left(a_{1},\ldots,a_{p}\right)\) where \(p\leq\frac{n}{2}\), we can find a half-diagram on \(n\) points which corresponds to this restricted sequence.  We define \(l_{i}\left(\mathfrak{a}\right)\) after Genauer and Stoltzfus \cite{2005math.....11003G} as the restricted sequence constructed from \(\mathfrak{a}\) by inserting an arch between the \(\left(i-1\right)\)th point and the \(i\)th point (where we consider the \(j\)th point of \(\mathfrak{a}\) to be the \(\left(j+2\right)\)th point of the new diagram, and the new arch connects the \(i\)th point and the \(\left(i+1\right)\)th point).

If the restricted sequence is empty, then it corresponds to the diagram on \(n\) points, each of whose points has a through-string.  Inductively, assume that any generalized restricted sequence with \(p\) integers for some \(p\leq\frac{n}{2}\) corresponds to a half-diagram on \(n\) points.  For a generalized restricted sequence \(\left(a_{1},\ldots,a_{p+1}\right)\), \(\left(a_{2},\ldots,a_{p+1}\right)\) is also a generalized restricted sequence, with \(p\) integers, which hence corresponds to a half-diagram \(\mathfrak{a}\) on \(n\) points.  Then \(l_{a_{1}}\left(\mathfrak{a}\right)\) is a half-diagram on \(n+2\) points.  The first closing bracket of \(\mathfrak{a}\) is the \(\left(a_{2}+1\right)\)th bracket, and \(a_{2}+1\geq a_{1}\), so if the \(j\)th bracket of \(\mathfrak{a}\) is a closing bracket, it becomes the \(\left(j+2\right)\)th bracket of \(l_{a_{1}}\left(\mathfrak{a}\right)\).  So any closing bracket of \(l_{a_{1}}\left(\mathfrak{a}\right)\), aside from the \(\left(a_{1}+1\right)\)th, has index greater than or equal to \(a_{1}+2\).  Thus the first closing bracket of \(l_{a_{1}}\left(\mathfrak{a}\right)\) is the \(\left(a_{1}\right)\)th, so its leftmost innnermost arch connects its \(a_{1}\)th point to its \(\left(a_{1}+1\right)\)th point, and hence the first integer in its restricted sequence is \(a_{1}\).  If we remove this arch, we get \(\mathfrak{a}\), which corresponds to the restricted sequence \(\left(a_{2},\ldots,a_{p}\right)\).  So the diagram \(l_{a_{1}}\left(\mathfrak{a}\right)\) corresponds to the restricted sequence \(\left(a_{1},\ldots,a_{p+1}\right)\).

Thus, we have a bijection between half-diagrams on \(n\) points and generalized restricted sequences with \(p\) points for \(p\leq\frac{n}{2}\).  We think of \(\mathfrak{a}\) being equal to the ordered \(p\)-tuple of integers which make up the corresponding restricted sequence.

We can then express the diagram corresponding to restricted sequence \(\left(a_{1}\ldots,a_{p}\right)\) as \(l_{a_{1}}\circ\ldots\circ l_{a_{p}}\left(\xi\right)\), where \(\xi\) is the diagram on \(n-p\) points where each point has a through-string.

At each stage of this construction, the closing bracket of the last arch inserted is the first closing bracket.  As the order of the brackets is not changed, the \(i\)th closing bracket belongs to the arch added by \(l_{a_{i}}\).  As the arch added by \(l_{a_{i}}\) is the leftmost innermost arch of \(l_{a_{i}}\circ\ldots\circ l_{a_{p}}\left(\xi\right)\), its closing bracket, the \(\left(a_{i}+1\right)\)th bracket, is the first closing bracket.  Then the height after the opening bracket of this arch or before its closing bracket is \(h_{a_{i}}\left(l_{a_{i}}\circ\ldots\circ l_{a_{p}}\left(\xi\right)\right)=a_{i}\).  We can interpret this as the base height inside this arch.  As the insertion of a new arch, that is, an opening bracket immediately followed by a closing bracket, cannot change the height at any of the other points, the height just after the opening bracket of the arch added by \(l_{a_{i}}\) or just before its closing bracket is still \(a_{i}\) (so the base height inside this arch does not change).  So we can also interpret the restricted sequence \(\left(a_{1},\ldots,a_{p}\right)\) as the height before each closing bracket, in order, left to right.

\subsection{Orderings on the Diagrams}

We define a partial order \(\preceq\) on half-diagrams with the same number of through-strings by letting \(\left(a_{1},\ldots,a_{p}\right)\preceq\left(b_{1},\ldots,b_{p}\right)\) if \(a_{i}\leq b_{i}\) for all \(i\), \(1\leq i\leq p\).  We define a linear order \(\leq\) extending this partial order by letting \(\left(a_{1}\ldots,a_{p}\right)\leq\left(b_{1},\ldots,b_{n}\right)\) if \(\left(a_{1}\ldots,a_{p}\right)=\left(b_{1},\ldots,b_{p}\right)\) or if there is an \(i\), \(1\leq i\leq n\) such that \(a_{i}<b_{i}\) and \(a_{j}=b_{j}\) for all \(j<i\) (lexicographical order on sequences of the same length).

In the context of Dyck paths, the partial order corresponds to {\em geometrical inclusion}; that is, \(\mathfrak{a}\preceq\mathfrak{b}\) iff the path corresponding to \(\mathfrak{a}\) is never above the path corresponding to \(\mathfrak{b}\).  More formally, \(\mathfrak{a}\preceq\mathfrak{b}\) if \(h_{i}\left(\mathfrak{a}\right)\leq h_{i}\left(\mathfrak{b}\right)\) for all \(i\), \(1\leq i\leq n\).

If there is a \(i\) such that \(h_{i}\left(\mathfrak{a}\right)>h_{i}\left(\mathfrak{b}\right)\), then we can pick the first such \(i\).  Then \(h_{i-1}\left(\mathfrak{a}\right)=h_{i-1}\left(\mathfrak{b}\right)\), so the \(i\)th step of \(\mathfrak{a}\) must be an upward step, while the \(i\)th step of \(\mathfrak{b}\) must be a downward step.  Let \(\mathfrak{a}=\left(a_{1},\ldots,a_{p}\right)\) and \(\mathfrak{b}=\left(b_{1},\ldots,b_{p}\right)\), and let this be the \(j\)th downward step of \(\mathfrak{b}\) (so there are \(j-1\) downward steps among the first \(i-1\) steps of \(\mathfrak{b}\), and hence of \(\mathfrak{a}\) as well).  Then \(b_{j}=h_{i-1}\left(\mathfrak{b}\right)=h_{i-1}\left(\mathfrak{a}\right)\).  \(a_{j}\) is equal to the height of \(\mathfrak{a}\) just before its \(j\)th closing bracket, which is its first closing bracket after the \(i\)th step.  As there may only be opening brackets between the \(i\)th step and its \(j\)th closing bracket, its height before its \(j\)th closing bracket must be greater than or equal to its height at \(i\) \(h_{i}\left(\mathfrak{a}\right)\), which is greater than \(h_{i}\left(\mathfrak{b}\right)\).  So \(a_{j}>b_{j}\).  So if \(a_{j}\leq b_{j}\) for all \(j\), \(1\leq j\leq p\), then \(h_{i}\left(\mathfrak{a}\right)\leq h_{i}\left(\mathfrak{b}\right)\) for all \(i\), \(1,\leq i\leq n\).

Conversely, if \(a_{j}>b_{j}\) for some \(j\), we can pick the first such \(j\).  There are \(j-1\) closing brackets before the \(j\)th closing bracket, so if the height right before the \(j\)th closing bracket of \(\mathfrak{a}\) is greater than the height right before the \(j\)th closing bracket of \(\mathfrak{b}\), then there must be more opening brackets in \(\mathfrak{a}\) before the \(j\)th bracket than in \(\mathfrak{b}\).  Thus the index of the \(j\)th bracket of \(\mathfrak{a}\) is greater than the index of the \(j\)th bracket of \(\mathfrak{b}\).  If the index of the \(j\)th bracket of \(\mathfrak{b}\) is \(i\), then there is one more closing bracket among the first \(i\) brackets of \(\mathfrak{b}\) than in \(\mathfrak{a}\), so \(h_{i}\left(\mathfrak{a}\right)>h_{i}\left(\mathfrak{b}\right)\).  So if \(h_{i}\left(\mathfrak{a}\right)\leq h_{i}\left(\mathfrak{b}\right)\) for all \(i\), \(1\leq i\leq n\), then \(a_{j}\leq b_{j}\) for all \(j\), \(1\leq j\leq p\).

\subsection{Enumerating Diagrams}

There are \(\binom{n}{p}\) sequences of \(n\) brackets with \(n-p\) opening brackets and \(p\) closing brackets, although these are not all legitimate bracket diagrams.  This is equivalent to relaxing the condition that the points in a Dyck path must have nonnegative \(y\)-coordinate.

\begin{figure}
\centering
\begin{multline*}
\includegraphics{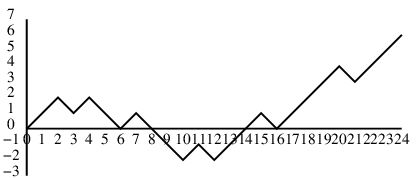}\mapsto\includegraphics{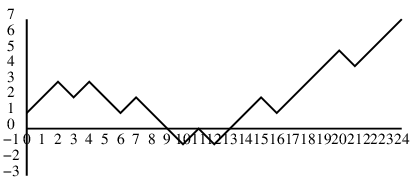}
\\\leftrightarrow\includegraphics{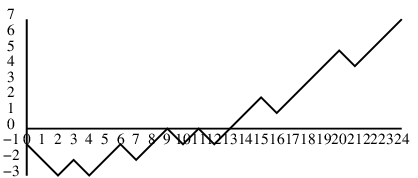}
\end{multline*}
\caption{The bijection between sequences between \(\left(0,0\right)\) and \(\left(n,n-2p\right)\) which are not legitimate and sequences between \(\left(0,-1\right)\) and \(\left(n,n-2p+1\right)\)}
\label{bijection}
\end{figure}

We can construct a bijection between the sequences which are not legitimate bracket diagrams and a similar set of sequences between \(\left(0,-1\right)\) and \(\left(n,n-2p+1\right)\) which we describe below and which we illustrate in Figure \ref{bijection}.

Given a sequence of points between \(\left(0,0\right)\) and \(\left(n,n-2p\right)\), we can add \(1\) to each of the \(y\)-coordinates, and get a sequence which begins at \(\left(0,1\right)\) and ends at \(\left(n,n-2p+1\right)\).  A sequence which had at least one point with a negative \(y\)-coordinate is now a sequence with at least one point whose \(y\)-coordinate is not positive.  As the value of the \(y\)-coordinate can only change by one, it must have at least one point whose \(y\)-coordinate is \(0\), and hence a first point, say the \(i\)th, with \(y\)-coordinate \(0\).

We can then construct a sequence from this sequence by changing the sign on the \(y\)-coordinate of all the points before the \(i\)th.  Then the first point is at \(\left(0,-1\right)\), since the \(y\)-coordinate of the first point is \(1\), and hence it is before the first point whose \(y\)-coordinate is \(0\).  Thus, this process uniquely defines a sequence of points from \(\left(0,-1\right)\) to \(\left(n,n-2p+1\right)\) where \(\left(x_{j+1},y_{j+1}\right)=\left(x_{j}+1,y_{j}\pm 1\right)\).

Conversely, given a sequence of points from \(\left(0,-1\right)\) to \(\left(n,n-2p+1\right)\) in which \(\left(x_{i+1},y_{i+1}\right)=\left(x_{i}+1,y_{i}\pm 1\right)\), there must be a point whose \(y\)-coordinate is \(0\), since there are points with negative \(y\)-coordinate and positive \(y\)-coordinate, and the \(y\)-coordinate may only change by \(\pm 1\) at each step.  So there must be a first such point.  We can perform the same process, changing the sign on the \(y\)-coordinate of each point before this one.  Then we have uniquely defined a sequence from \(\left(0,1\right)\) to \(\left(n,n-2p+1\right)\) where \(\left(x_{i+1},y_{i+1}\right)=\left(x_{i}+1,y_{i}\pm 1\right)\), and the first point with \(0\) \(y\)-coordinate is the same point.  Thus, performing the above process on this new path gives us the original path back again.

Thus, we have a bijection between the paths from \(\left(0,1\right)\) to \(\left(n,n-2p+1\right)\) which have a point with \(0\) \(y\)-coordinate and those from \(\left(0,-1\right)\) to \(\left(n,n-2p\right.\allowbreak\left.+1\right)\).  We know that there are \(\binom{n}{p-1}\) of the latter.  The former are also in bijection with sequences of brackets which are not legitimate.  Thus, there are \(\binom{n}{p}-\binom{n}{p-1}=\binom{n-1}{p}-\binom{n-1}{p-2}\) legitimate generalized Dyck paths from \(\left(0,0\right)\) to \(\left(n,n-2p\right)\).

After \cite{MR1608551}, we define \(c_{n,n-2p}=\binom{n}{p}-\binom{n}{2p-1}=\binom{n-1}{p}-\binom{n-1}{p-2}\), and call it a generalized Catalan number.

\section{Chebyshev Polynomials}

We define the \(n\)th Chebyshev polynomial \(\Delta_{n}\left(q\right)\) as the determinant of the \(n\times n\) matrix whose diagonal entries are all \(q\) and whose subdiagonal entries and superdiagonal entries are all \(1\).  So
\[\Delta_{n}\left(q\right)=\left|\begin{array}{ccccc}q&1&0&\cdots&0\\1&q&\ddots&\ddots&\vdots\\0&\ddots&\ddots&\ddots&0\\\vdots&\ddots&\ddots&q&1\\0&\cdots&0&1&q\end{array}\right|\textrm{.}\]
The Chebyshev polynomials satisfy the recurrence relation
\[\Delta_{n+1}\left(q\right)=q\Delta_{n}\left(q\right)-\Delta_{n-1}\left(q\right)\textrm{.}\]
If we define \(\Delta_{-1}\left(0\right)\) and \(\Delta_{0}\left(q\right)=1\), these also satisfy the recurrence relation.  The first few Chebyshev polynomials are \[\begin{array}{l}\Delta_{-1}\left(q\right)=0\\\Delta_{0}\left(q\right)=1\\\Delta_{1}\left(q\right)=q\\\Delta_{2}\left(q\right)=q^{2}-1\\\Delta_{3}\left(q\right)=q^{3}-2q\\\Delta_{4}\left(q\right)=q^{4}-3q^{2}+1\\\Delta_{5}\left(q\right)=q^{5}-4q^{3}+3q\end{array}\]

We define \(\mu_{i}=\frac{\Delta_{i-1}}{\Delta_{i}}\).  In this paper, we only deal with values of \(q\) for which all relevant Chebyshev polynomials are not zero, so this quotient will always be well defined.

\section{The Orthogonal Basis}

We define an orthogonal basis, following the method of box addition in \cite{MR1608551} and \cite{MR1999738}.  However, we do not normalize our basis, so our definition of box addition differs from the definitions in \cite{MR1608551} and \cite{MR1999738} by a scalar factor.  \cite{MR1462755} presents both the normalized and the unnormalized bases.

\subsection{Box Addition}

\begin{figure}
\centering
\includegraphics{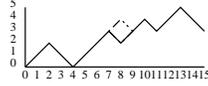}
\caption{Box addition at \(8\) performed on a Dyck path}
\label{box}
\end{figure}

We define {\em box addition at position \(i\)}, denoted \(\Diamond_{i}\), on Dyck path \(\mathfrak{a}\) which has a minimum at \(i\).  We let \(\Diamond_{i}\left(\mathfrak{a}\right)\) be the diagram constructed from \(\mathfrak{a}\) by replacing the downward step at \(i\) with an upward step and the upward step at \(i+1\) with a downward step, as in Figure \ref{box}.  \(\Diamond_{i}\) is not defined when \(\mathfrak{a}\) does not have a minimum at \(i\).

\(\mathfrak{a}\) and \(\Diamond_{i}\left(\mathfrak{a}\right)\) have the same sequence of opening and closing brackets up to the \(i\)th bracket, so for \(j<i\), \(h_{j}\left(\Diamond_{i}\left(\mathfrak{a}\right)\right)=h_{j}\left(\mathfrak{a}\right)\).  As an opening bracket has been replaced by a closing bracket and a closing bracket by an opening bracket in the construction of \(\Diamond_{i}\left(\mathfrak{a}\right)\), the number of opening brackets and closing brackets in the first \(j\) brackets of \(\mathfrak{a}\) and \(\Diamond_{i}\left(\mathfrak{a}\right)\) are equal for \(j>i\).  So for \(j>i\), \(h_{j}\left(\Diamond_{i}\left(\mathfrak{a}\right)\right)=h_{j}\left(\mathfrak{a}\right)\).  Thus \(h_{j}\left(\Diamond_{i}\left(\mathfrak{a}\right)\right)=h_{j}\left(\mathfrak{a}\right)\) for all \(j\neq i\).  As \(\Diamond_{i}\left(\mathfrak{a}\right)\) has one fewer closing bracket and one more opening bracket in its first \(i\) brackets, \(h_{i}\left(\Diamond_{i}\left(\mathfrak{a}\right)\right)=h_{i}\left(\mathfrak{a}\right)+2\).  So \(h_{j}\left(\Diamond_{i}\left(\mathfrak{a}\right)\right)\geq h_{j}\left(\mathfrak{a}\right)\) for all \(j\), \(1\leq j\leq n\).  (This also shows that box addition respects the orderings, so if \(\mathfrak{a}\preceq\mathfrak{b}\), then \(\Diamond_{i}\left(\mathfrak{a}\right)\preceq\Diamond_{i}\left(\mathfrak{b}\right)\), if both are well-defined.)

In fact, if \(\mathfrak{a}=\left(a_{1},\ldots,a_{p}\right)\), then the height of \(\mathfrak{a}\) just before its \(j\)th closing bracket is \(a_{j}\).  If the closing bracket at \(i\) is the \(j\)th closing bracket of \(\mathfrak{a}\), then the closing bracket at \(i+1\) in \(\Diamond_{i}\left(\mathfrak{a}\right)\) is the \(j\)th closing bracket of \(\Diamond_{i}\left(\mathfrak{a}\right)\).  As there are \(j-1\) closing brackets before each, and one more opening bracket before the \(j\)th closing bracket of \(\mathfrak{b}\) than before the \(j\)th closing bracket of \(\mathfrak{a}\).  Thus \(b_{j}=a_{j}+1\).  As the rest of the closing brackets of the diagrams are at the same indices, and none aside from the \(j\)th closing bracket of \(\Diamond_{i}\left(\mathfrak{a}\right)\) is preceded by the \(i\)th point, the height of the diagrams right before each other closing bracket in each diagram is the same, and hence so are the integers in their restricted sequences.  So \(\Diamond_{i}\left(\mathfrak{a}\right)\) corresponds to restricted sequence \(\left(a_{1},\ldots,a_{j-1},a_{j}+1,a_{j+1},\ldots,a_{p}\right)\).

\begin{figure}
\centering
\includegraphics{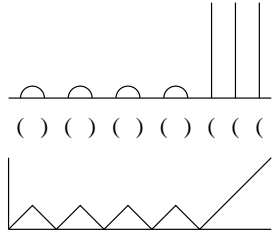}
\caption{The minimal element with four pairs and three through-strings, \(\left(1,1,1,1\right)\), shown as a half-pairing, a bracket diagram and a Dyck path}
\label{minimal}
\end{figure}

The \(p\)-tuple of positive integers \(\left(1,\ldots,1\right)\) is a legitimate restricted sequence.  Each index is less than or equal to the corresponding index in any other sequence of positive integers, so this is smaller than any other restricted sequence in the partial order.  We call this the {\em minimal} or {\em fundamental} element.  Figure \ref{minimal} shows the half-diagram, generalized Dyck path and bracket diagram corresponding to \(\left(1,\ldots,1\right)\).

If \(\left(a_{1},\ldots,a_{p}\right)\) and \(\left(b_{1},\ldots,b_{p}\right)\) are restricted sequences with \(\left(a_{1},\ldots,a_{p}\right)\prec\left(b_{1},\ldots,b_{p}\right)\), then there must be at least one \(i\) such that \(a_{i}<b_{i}\).  Let \(i\) be the largest such index.  Then we can let \(c_{i}=a_{i}+1\leq b_{i}\) and \(c_{j}=a_{j}\leq b_{j}\) for all \(j\neq i\), \(1\leq j\leq p\).  Then \(c_{i}=a_{i}+1\geq a_{i-1}+1>a_{i-1}=c_{i-1}\) if \(i-1\) is a valid index, \(c_{i+1}=a_{i+1}=b_{i+1}\geq b_{i}-1\geq c_{i}-1\) if \(i+1\) is a valid index, and \(c_{j+1}=a_{j+1}\geq a_{j}-1=c_{j}-1\) for all \(j\) not equal to \(i-1\) or \(i\).  Thus \(\left(c_{1},\ldots,c_{p}\right)\) is a valid restricted sequence, \(\left(b_{1},\ldots,b_{p}\right)\) is the result of box addition on \(\left(c_{1},\ldots,c_{p}\right)\), and \(\left(a_{1},\ldots,a_{p}\right)\prec\left(c_{1},\ldots,c_{p}\right)\preceq\left(b_{1},\ldots,b_{p}\right)\).  Inductively, we can see that we can construct \(\left(b_{1},\ldots,b_{p}\right)\) from \(\left(a_{1},\ldots,a_{p}\right)\) by repeated box additions.

In particular, all diagrams can be constructed from \(\left(1,\ldots,1\right)\) by a sequence of box additions.  Furthermore, the number of box additions required is always \(a_{1}+\ldots+a_{p}-p\) (since this quantity is zero on \(\left(1,\ldots,1\right)\), and increases by one with each box addition), regardless of which box additions we perform, and in what order, so it is an invariant of the diagram.  This allows us to perform induction on this quantity, the {\em number of boxes} of \(\mathfrak{a}\), or \(\#_{\Diamond}\left(\mathfrak{a}\right)\) ({\em length} or \(\left|\mathfrak{a}\right|\) in \cite{MR1608551}), which we do in several of the proofs.

\subsubsection{Grey Box Addition}

\begin{figure}
\centering
\[\includegraphics{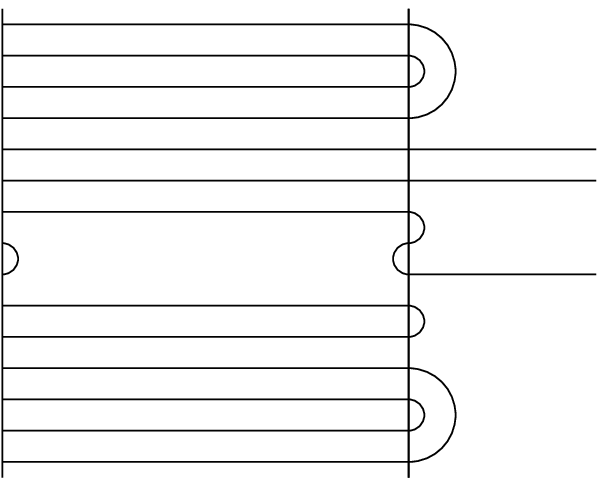}=\includegraphics{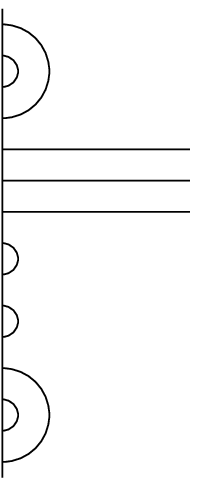}\]
\caption{The action of \(e_{i}\) on \(\xi_{\mathfrak{a}}\) maps it to \(\xi_{\Diamond_{i}\left(\mathfrak{a}\right)}\), the grey box addition equivalent to the box addition shown in Figure \ref{box}}
\label{greybox}
\end{figure}

\begin{prop}
If \(\mathfrak{a}\) has a minimum at \(i\), then \(\xi_{\Diamond_{i}\left(\mathfrak{a}\right)}=e_{i}\xi_{\mathfrak{a}}\).

\begin{proof}
We show this diagrammatically in Figure \ref{greybox}.

We know that \(e_{i}\xi_{\mathfrak{a}}=c\xi_{\mathfrak{b}}\) for some diagram \(\mathfrak{b}\) and some scalar \(c\), which may be \(0\) or a power of \(q\).

The \(j\)th point of \(e_{i}\), for \(j\) not equal to \(i\) or \(i+1\), is connected to the \(\left(n-j+1\right)\)th point of \(e_{i}\), which is identified with the \(j\)th point of \(\xi_{\mathfrak{a}}\).  If the \(j\)th point of \(\mathfrak{a}\) was not paired with the \(i\)th point or the \(\left(i+1\right)\)th point, then the \(j\)th point of \(\mathfrak{b}\) must be paired with the point that the \(j\)th point of \(\xi_{\mathfrak{a}}\) was paired with.  Thus any step of the resulting diagram is an upward or downward step if the corresponding step in \(\mathfrak{a}\) is an upward or downward step.

The \(i\)th point of \(e_{i}\) is connected to its \(\left(i+1\right)\)th point, so the \(i\)th point of \(e_{i}\xi_{\mathfrak{a}}\) is connected to its \(\left(i+1\right)\)th point.  Then the \(i\)th point is connected to a point of greater index and the \(\left(i+1\right)\)th point is connected to a point of smaller index, so the \(i\)th step is an upward step and the \(\left(i+1\right)\) is a downward step, as they are in \(\Diamond_{i}\left(\mathfrak{a}\right)\).

The only remaining points are those that were originally paired with \(i\) and \(i+1\) in \(\mathfrak{a}\), so these must be paired with each other.  As the \(i\)th step of \(\mathfrak{a}\) is decreasing, its \(i\)th point is paired with a point of smaller index, and since the \(\left(i+1\right)\)th step of \(\mathfrak{a}\) is increasing, it is paired with a point of greater index.  So the original partner of \(i\), whose index is less than \(i\) and hence less than \(i+1\) and hence less than the index of the partner of \(i+1\) in \(\mathfrak{a}\), is paired in \(\mathfrak{b}\) with a point of greater index.  So it corresponds to an increasing step, and the point partnered with \(i+1\) in \(\mathfrak{a}\) corresponds to a decreasing step.  As the \(i\)th point of \(\mathfrak{a}\) was joined in \(\mathfrak{a}\) to a point of smaller index (and the \(\left(i+1\right)\)th point was joined to a point of greater index), its partner corresponds to an increasing step (and the partner of the \(\left(i+1\right)\)th point corresponds to a decreasing step).  So these steps in \(\mathfrak{b}\) are the same as those in \(\mathfrak{a}\).

Any sequence of edges joining the point at infinity back to itself or any closed loop which occurs when we have identified the \(\left(n-j+1\right)\)th point of \(e_{i}\) with the \(j\)th point of \(\xi_{\mathfrak{a}}\) must involve at least one arch of \(e_{i}\), since the arch diagram corresponding to \(\mathfrak{a}\) contains no such structures.  This arch must connect points in \(\left\{n+1,\ldots,2n\right\}\), since the points in \(\left\{1,\ldots,n\right\}\) have only one edge.  There is only one such arch in \(e_{i}\), namely, that connecting the \(\left(2n-i+1\right)\)th point to the \(\left(2n-i\right)\)th.  These points are identified with the \(i\)th and \(\left(i+1\right)\)th points of \(\xi_{\mathfrak{a}}\).  Since \(\mathfrak{a}\) has a minimum at \(i\), the \(i\)th point is connected to a point of smaller index, and the \(\left(i+1\right)\)th to a point of greater index, they are both connected to points which are not \(i\) or \(i+1\) in \(\mathfrak{a}\).  These are identified with points of the form \(2n-j+1\), \(1\leq j\leq n\) in \(e_{i}\), for \(j\) not equal to \(i\) or \(i+1\).  These points are connected in \(e_{i}\) to the \(j\)th point, so these arches cannot be part of any sequence of edges connecting the point at infinity back to itself or any closed loop.  So there is no such structure.  So \(c\) cannot be equal to zero, and since there are no closed loops, \(c=q^{0}=1\).  So \(e_{i}\xi_{\mathfrak{a}}=\xi_{\Diamond_{i}\left(\mathfrak{a}\right)}\).
\end{proof}
\end{prop}

We call this construction grey box addition, after Di Francesco (\cite{MR1608551}).  This allows us to construct our original basis by letting \(TL_{n}\) act on the basis element of \(U\left(n;p\right)\) corresponding to the minimal diagram \(\xi_{\left(1,\ldots,1\right)}\), which shows us that each \(U\left(n;p\right)\) is not reducible.

\subsubsection{White Box Addition}

We define our second basis in terms of a similar operation, which we will call white box addition, again after Di Francesco (\cite{MR1608551}).  This basis is also indexed by diagrams.  We will denote members of the second basis with a prime, \(\xi^{\prime}_{\mathfrak{a}}\).  We define \(\xi^{\prime}_{\left(1,\ldots,1\right)}=\xi_{\left(1,\ldots,1\right)}\).  If \(\mathfrak{a}\) has a minimum at \(i\), then we define \(\xi^{\prime}_{\Diamond_{i}\left(\mathfrak{a}\right)}=\left(e_{i}-\mu_{h_{i}\left(\mathfrak{a}\right)+1}\right)\xi^{\prime}_{\mathfrak{a}}\).

\begin{prop}
\(\xi^{\prime}_{\mathfrak{a}}\) is a linear combination of \(\xi_{\mathfrak{b}}\), where \(\mathfrak{b}\preceq\mathfrak{a}\).
\begin{proof}
\(\xi^{\prime}_{\left(1,\ldots,1\right)}=\xi_{\left(1,\ldots,1\right)}\), so our proposition holds for the minimal element.

Assume that our proposition holds for some \(\mathfrak{a}\); that is, \(\xi^{\prime}_{\mathfrak{a}}\) is a linear combination of \(\xi_{\mathfrak{b}}\) for \(\mathfrak{b}\preceq\mathfrak{a}\).  Then \(\xi^{\prime}_{\Diamond_{i}\left(\mathfrak{a}\right)}=\left(e_{i}-\mu_{h_{i}\left(\mathfrak{a}\right)+1}\right)\xi^{\prime}_{\mathfrak{a}}=e_{i}\xi^{\prime}_{\mathfrak{a}}-\mu_{h_{i}\left(\mathfrak{a}\right)+1}\xi^{\prime}_{\mathfrak{a}}\).

The latter term is a linear combination of \(\xi_{\mathfrak{b}}\) where \(\mathfrak{b}\preceq\mathfrak{a}\).  We examine the action of \(e_{i}\) on such a \(\xi_{\mathfrak{b}}\).

We have seen above that if \(\mathfrak{b}\) has a minimum at \(i\), \(e_{i}\xi_{\mathfrak{b}}=\xi_{\Diamond_{i}\left(\mathfrak{b}\right)}\).  If \(\mathfrak{b}\preceq\mathfrak{a}\), then \(\Diamond_{i}\left(\mathfrak{b}\right)\preceq\Diamond_{i}\left(\mathfrak{a}\right)\).

If \(\mathfrak{b}\) has a maximum at \(i\), then it has an opening bracket at \(i\) and a closing bracket at \(i+1\), and hence an arch connecting its \(i\)th point to its \(\left(i+1\right)\)th.  These points are identified with the \(\left(2n-i+1\right)\)th and \(\left(2n-i\right)\)th of \(e_{i}\), which are themselves connected.  So this forms a closed loop, and we get a factor of \(q\).  Since the \(i\)th point of \(e_{i}\) is connected to its \(\left(i+1\right)\)th point, the \(i\)th point of the new diagram is connected to its \(\left(i+1\right)\)th point, as in \(\xi_{\mathfrak{b}}\).  If \(j\) is not equal to \(i\) or \(i+1\), \(1\leq j\leq n\), then the \(j\)th point of \(e_{i}\) is connected to its \(\left(2n-j+1\right)\)th point, which is identified with the \(j\)th point of \(\xi_{\mathfrak{b}}\).  This point is connected to some other point which is also not equal to \(i\) or \(i+1\), since these are connected to each other.  So the \(j\)th point of \(e_{i}\xi_{\mathfrak{b}}\) is connected to the point of the same index as the partner of the \(j\)th point of \(\xi_{\mathfrak{b}}\).  So \(e_{i}\xi_{\mathfrak{b}}\) has all the same pairing as \(\mathfrak{b}\), and one closed loop.  As the one loop connecting two points in \(\left\{n+1,\ldots,2n\right\}\) in \(e_{i}\) is part of the closed loop, there are no other closed loops or sequences of edges connecting infinity back to itself.  So \(e_{i}\xi_{\mathfrak{b}}=q\xi_{\mathfrak{b}}\), and \(\mathfrak{b}\preceq\mathfrak{a}\preceq\Diamond_{i}\left(\mathfrak{a}\right)\).

If \(\mathfrak{b}\) has a slope at \(i\), then it has two brackets of the same type at \(i\) and \(i+1\).  \(e_{i}\xi_{\mathfrak{b}}=c\xi_{\mathfrak{c}}\) where \(c\) is either \(0\) or a power of \(q\), and \(\mathfrak{c}\) is a half-diagram.  If the two brackets corresponding to through-strings, then the arch from \(2n-i+1\) to \(2n-i\) in \(e_{i}\) connects them, so \(c=0\).  Otherwise, \(\mathfrak{c}\) has an opening bracket at \(i\) and a closing bracket at \(i+1\), where \(\mathfrak{b}\) had either two opening brackets or two closing brackets.  \(\mathfrak{c}\) must have the same number of opening brackets and closing brackets as \(\mathfrak{b}\), and the points which were not either \(i\) or \(i+1\) or partnered with one of these are paired just as they are in \(\mathfrak{b}\), so the type of bracket at any of these points is the same as in \(\mathfrak{b}\).  The type of bracket at one of \(i\) and \(i+1\) was changed, so the type of bracket at one of their partners (which would be of opposite type) must be changed as well.  If the brackets at \(i\) and \(i+1\) are both opening brackets, then their partners are of greater index, and if the brackets at \(i\) and \(i+1\) are both closing brackets, then their partners must be of smaller index.  Either way, we can construct \(\mathfrak{c}\) from \(\mathfrak{b}\) by changing an opening bracket to a closing bracket at one index, and a closing bracket to an opening bracket at a greater index.  So the height of \(\mathfrak{c}\) at any index must be less than or equal to the height of \(\mathfrak{b}\) at that index.  So \(\mathfrak{c}\preceq\mathfrak{b}\preceq\mathfrak{a}\).

Since \(\xi^{\prime}_{\mathfrak{a}}\) is a linear combination of \(\xi_{\mathfrak{b}}\) with \(\mathfrak{b}\preceq\mathfrak{a}\) and \(e_{i}\xi_{\mathfrak{b}}\) produces a linear combination of \(\xi_{\mathfrak{c}}\) with \(\mathfrak{c}\preceq\Diamond_{i}\left(\mathfrak{a}\right)\), \(e_{i}\xi^{\prime}_{\mathfrak{a}}\) is a linear combination of \(\xi_{\mathfrak{b}}\) with \(\mathfrak{b}\preceq\mathfrak{a}\), and hence so is \(\left(e_{i}-\mu_{h_{i}\left(\mathfrak{a}\right)+1}\right)\xi^{\prime}_{\mathfrak{a}}\).
\end{proof}
\end{prop}

\begin{prop}
The coefficient of \(\xi_{\mathfrak{a}}\) in \(\xi^{\prime}_{\mathfrak{a}}\) is \(1\).
\begin{proof}
Clearly, the coefficient of \(\xi_{\left(1,\ldots,1\right)}\) in \(\xi^{\prime}_{\left(1,\ldots,1\right)}=\xi_{\left(1,\ldots,1\right)}\) is \(1\).

Assume that the coefficient of \(\xi_{\mathfrak{a}}\) in \(\xi^{\prime}_{\mathfrak{a}}\) is also \(1\).  By definition, \(\xi^{\prime}_{\Diamond_{i}\left(\mathfrak{a}\right)}=\left(e_{i}-\mu_{h_{i}\left(\mathfrak{a}\right)+1}\right)\xi^{\prime}_{\mathfrak{a}}=e_{i}\xi^{\prime}_{\mathfrak{a}}-\mu_{h_{i}\left(\mathfrak{a}\right)+1}\xi^{\prime}_{\mathfrak{a}}\).  The second term is a linear combination of \(\xi_{\mathfrak{b}}\) with \(\mathfrak{b}\preceq\mathfrak{a}\prec\Diamond_{i}\left(\mathfrak{a}\right)\), so this does not contribute to the coefficient on \(\xi_{\mathfrak{a}}\).  The coefficient of \(\xi_{\mathfrak{a}}\) in \(\xi^{\prime}_{\mathfrak{a}}\) is \(1\), so the term \(e_{i}\xi_{\mathfrak{a}}=\xi_{\Diamond_{i}\left(\mathfrak{a}\right)}\) appears in \(e_{i}\xi^{\prime}_{\mathfrak{a}}\).

If \(\mathfrak{b}\prec\mathfrak{a}\) has a minimum at \(i\), then \(e_{i}\xi_{\mathfrak{b}}=\xi_{\Diamond_{i}\left(\mathfrak{b}\right)}\), and since \(\mathfrak{b}\prec\mathfrak{a}\), \(\Diamond_{i}\left(\mathfrak{b}\right)\prec\Diamond_{i}\left(\mathfrak{b}\right)\).  If \(\mathfrak{b}\prec\mathfrak{a}\) has a maximum at \(i\), then \(e_{i}\xi_{\mathfrak{b}}=q\xi_{\mathfrak{b}}\).  If \(\mathfrak{b}\prec\mathfrak{a}\) has a slope at \(i\), then \(e_{i}\xi_{\mathfrak{b}}=\xi_{\mathfrak{c}}\), where \(\mathfrak{c}\preceq\mathfrak{b}\prec\mathfrak{a}\).

In any of these cases, \(e_{i}\xi_{\mathfrak{b}}\) is a linear combination of \(\xi_{\mathfrak{c}}\) with \(\mathfrak{c}\prec\mathfrak{a}\), so such terms do not contribute to the coefficient of \(\xi_{\mathfrak{a}}\) in \(\xi^{\prime}_{\mathfrak{a}}\).  So this coefficient is \(1\).
\end{proof}
\end{prop}

Thus the transformation \(\xi_{\mathfrak{a}}\mapsto\xi^{\prime}_{\mathfrak{a}}\) is an upper triangular transformation (with respect to the partial order, and hence also with respect to the linear order which extends it) with \(1\)s on the diagonal.  So it is invertible, and hence a legitimate change of basis.  So the \(\xi^{\prime}_{\mathfrak{a}}\) form a basis for \(U\left(n;p\right)\).

\subsubsection{Commutativity Relations}

We now demonstrate that white box addition is well-defined; that is, we get consistent definitions of \(\xi^{\prime}_{\mathfrak{a}}\) regardless of the order in which the box addition is done.  Grey box addition produces the original basis elements corresponding to appropriate Dyck paths, so it must be well-defined.  In either case, if \(\Diamond_{i}\left(\Diamond_{j}\left(\mathfrak{a}\right)\right)=\Diamond_{j}\left(\Diamond_{i}\left(\mathfrak{a}\right)\right)\) for \(i\neq j\), then \(\mathfrak{a}\) must have minima at \(i\) and \(j\).  Both the \(i\)th and \(j\)th steps must then be downward steps and the \(\left(i+1\right)\)th and \(\left(j+1\right)\)th steps must be upward steps, so \(j\neq i\pm 1\), or \(\left|i-j\right|>1\).  According to the Jones relations, \(e_{i}\) and \(e_{j}\) commute, and hence the factors corresponding to either type of box addition also commute.  So the box addition operations satisfy the same commutativity relations as box additions on a Dyck path.  We show more formally that our definition of the second basis of the \(\xi^{\prime}_{\mathfrak{a}}\) is well-defined.

\begin{prop}
If \(\mathfrak{b}\) can be derived from \(\mathfrak{a}\) by two distinct sequences of box additions, then the expressions for \(\xi^{\prime}_{\mathfrak{b}}\) corresponding to these two sequences are equal.
\begin{proof}
As shown above, if \(\mathfrak{b}\) can be derived from \(\mathfrak{a}\) by two distinct sequences of box additions, then the number of box additions is the same in each case, as shown above.  Furthermore, there must be the same number of box additions at each \(i\), \(1,\leq i\leq n\), since each box addition at \(i\) increases the height at \(i\) by exactly two, and no box addition at any other index alters the height at \(i\) (so there must be exactly \(\left(h_{i}\left(\mathfrak{b}\right)-h_{i}\left(\mathfrak{a}\right)\right)/2\) box additions at \(i\)).  Finally, the heights at which the box additions at \(i\) occur are the same; that is, there must be one at height \(h_{i}\left(\mathfrak{a}\right)\), one at \(h_{i}\left(\mathfrak{a}\right)+2\), etc., and they must occur in this order.  The white box addition operator corresponds to a factor of \(e_{i}-\mu_{h_{i}\left(\mathfrak{a}\right)+1}\), multiplied on the left, which depends only on the index \(i\) at which the box addition takes place and the height of the Dyck path at \(i\) when the box addition is performed, so the same factors must appear in our two expressions for \(\xi^{\prime}_{\mathfrak{b}}\), although possibly in a different order.

We prove the theorem by induction on the number of box additions, and hence the number of factors.  If there is only one factor in each sequence, then as established above, it must be the same factor, so the expressions are equal.

We now assume that, if there are \(k\) box additions in each of the two sequences, and hence \(k\) factors in each expression, then the two expressions are equal.  Throughout, we consider the first factor to be the first applied to the half-diagram, that is, the rightmost, rather than the first that would appear in a written expression of the product.  Let \(\mathfrak{b}\) be the result of \(k+1\) box additions on \(\mathfrak{a}\) in two possibly different ways.  Let the first box addition in the first sequence be \(\Diamond_{i}\).  Then \(\mathfrak{b}\) can be expressed as \(k\) box additions on \(\Diamond_{i}\left(\mathfrak{a}\right)\).  \(\Diamond_{i}\) corresponds to a factor of \(e_{i}-\mu_{h_{i}\left(\mathfrak{a}\right)+1}\).  We know that this factor occurs in the second expression, before any other factors corresponding to any other box addition at \(i\).  \(\mathfrak{a}\) has a minimum at \(i\), so the \(i\)th step of \(\mathfrak{a}\) is a downward step, and its \(\left(i+1\right)\)th step is an upward step, and hence no box addition can be performed at \(i-1\) or \(i+1\).  As no other box addition aside from one at \(i\) will change the types of steps at \(i\) and \(i+1\), there may be no box additions at \(i-1\) or \(i+1\) before the first box addition at \(i\).  So all of the factors that precede \(e_{i}-\mu_{h_{i}\left(\mathfrak{a}\right)+1}\) are of the form \(e_{j}-\mu_{h}\) for some \(j\) with \(\left|i-j\right|>1\) and some \(h\).  By the Jones relations, \(e_{i}\) commutes with \(e_{j}\), and hence the factor corresponding to \(\Diamond_{i}\) commutes with all of the factors that precede it.  Thus we can move the factor \(e_{i}-\mu_{h_{i}\left(\mathfrak{a}\right)+1}\) to the beginning.  The other \(k\) factors correspond to the \(k\) remaining box additions, which, as discussed above, are the same box additions at the same heights as the \(k\) that are applied to \(\Diamond_{i}\left(\mathfrak{a}\right)\) to get \(\mathfrak{b}\).  These must be in a legitimate order to correspond to these box additions.  Those that originally occurred before the first factor of \(e_{i}-\mu_{h_{i}\left(\mathfrak{a}\right)}\) and which now occur after it did not occur at \(i-1\), \(i\) or \(i+1\), so the existence of a minimum at their index is not affected.  For those that originally occurred after the first factor of \(e_{i}-\mu_{h_{i}\left(\mathfrak{a}\right)}\), exactly the same set of box additions at the same heights have been performed, so they are still performed legitimately.  Thus the second expression for \(\xi^{\prime}_{\mathfrak{a}}\) is equal to the expression for \(\xi^{\prime}_{\mathfrak{b}}\) corresponding to \(k\) box additions applied to \(\Diamond_{i}\left(\mathfrak{a}\right)\).  As these are the same \(k\) box additions as in the first sequence, both expressions for \(\xi^{\prime}_{\mathfrak{b}}\) are expressions we would get by applying the same \(k\) box additions to \(\Diamond_{i}\left(\mathfrak{a}\right)\).  By our inductive hypothesis, these two expressions for \(\xi^{\prime}_{\mathfrak{b}}\) are equal.

Thus our definition of the second basis of the \(\xi^{\prime}_{\mathfrak{a}}\) is consistent, regardless of the order in which the box additions are performed.
\end{proof}
\end{prop}

\subsection{Uniqueness of the Inner Product}

\begin{thm}
Up to scalar multiples, the inner product we have defined is the only inner product on \(U\left(n\right)\) such that the involution \(*\) on \(TL_{n}\) is an adjoint.
\begin{proof}

Assume there is another sesquilinear function \(\langle\cdot,\cdot\rangle^{\prime}\) for which the involution \(*\) is an adjoint.

\begin{figure}
\centering
\[q^{-6}\includegraphics{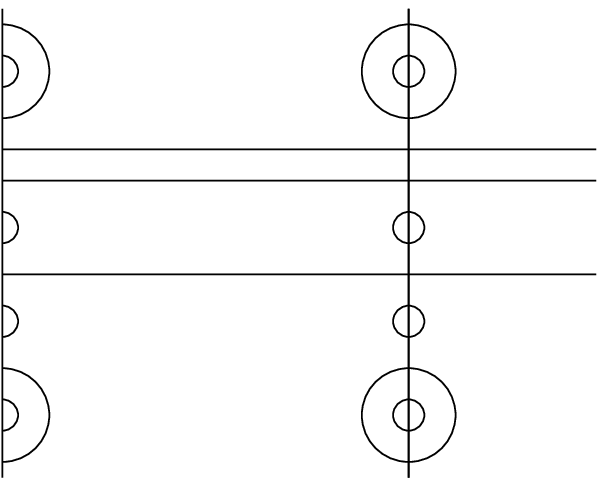}=\includegraphics{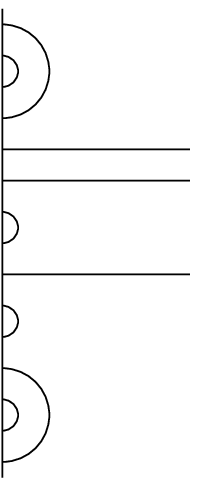}\]
\caption{We use \(e_{\mathfrak{a}\otimes\mathfrak{a}}\in TL_{n}\) to externalize the structure of \(\mathfrak{a}\)}
\label{externalization}
\end{figure}

For any arch connecting \(i\) and \(j\) in \(\xi_{\mathfrak{a}}\), the diagram \(e_{\mathfrak{a}\otimes\mathfrak{a}}\xi_{\mathfrak{a}}\) has an arch connecting \(i\) and \(j\), and for every through-string at \(i\) in \(\xi_{\mathfrak{a}}\), there is an arch connecting \(i\) and \(2n-i\) in \(e_{\mathfrak{a}\otimes\mathfrak{a}}\).  As the \(\left(2n-i\right)\)th point of \(e_{\mathfrak{a}\otimes\mathfrak{a}}\) is identified with the \(i\)th point of \(\xi_{\mathfrak{a}}\) in the diagram \(e_{\mathfrak{a}\otimes\mathfrak{a}}\xi_{\mathfrak{a}}\), the product has a through-string wherever \(\xi_{\mathfrak{a}}\) does.  Finally, for each arch of \(\xi_{\mathfrak{a}}\) connecting \(i\) and \(j\), there is an arch in \(e_{\mathfrak{a}\otimes\mathfrak{a}}\) connecting \(2n-i\) and \(2n-j\).  As the \(\left(2n-i\right)\)th point of \(e_{\mathfrak{a}\otimes\mathfrak{a}}\) is identified with the \(i\)th point of \(\xi_{\mathfrak{a}}\), both ends of these arches are identified and form a closed loop.  So for each of the \(p\) arches of \(\xi_{\mathfrak{a}}\), there is a closed loop in the diagram \(e_{\mathfrak{a}\otimes\mathfrak{a}}\xi_{\mathfrak{a}}\).  So \(e_{\mathfrak{a}\otimes\mathfrak{a}}\xi_{\mathfrak{a}}=q^{p}\xi_{\mathfrak{a}}\).  This method of externalizing the structure of \(\mathfrak{a}\) is shown in Figure \ref{externalization}.

Thus \(\xi_{\mathfrak{a}}=q^{-p}e_{\mathfrak{a}\otimes\mathfrak{a}}\xi_{\mathfrak{a}}\).  For each arch joining points \(i\) and \(j\) in \(e_{\mathfrak{a}\otimes\mathfrak{a}}\), there is an arch joining points \(2n-i\) and \(2n-j\), so \(e_{\mathfrak{a}\otimes\mathfrak{a}}^{*}=e_{\mathfrak{a}\otimes\mathfrak{a}}\), and hence \(\langle\xi_{\mathfrak{a}},\xi_{\mathfrak{b}}\rangle^{\prime}=\langle q^{-p}e_{\mathfrak{a}\otimes\mathfrak{a}}\xi_{\mathfrak{a}},\xi_{\mathfrak{b}}\rangle^{\prime}=q^{-p}\langle\xi_{\mathfrak{a}},e_{\mathfrak{a}\otimes\mathfrak{a}}\xi_{\mathfrak{b}}\rangle^{\prime}\).

In the diagram \(\langle\xi_{\mathfrak{a}},\xi_{\mathfrak{b}}\rangle\), the \(i\)th point of \(\xi_{\mathfrak{b}}\) is identified with the \(i\)th point of \(\xi_{\mathfrak{a}}\), and in the diagram \(e_{\mathfrak{a}\otimes\mathfrak{a}}\xi_{\mathfrak{a}}\), the \(i\)th point of \(\xi_{\mathfrak{b}}\) is identified with the \(\left(2n-i\right)\)th point of \(e_{\mathfrak{a}\otimes\mathfrak{a}}\).  For any arch connecting the \(i\)th and \(j\)th points of \(\xi_{\mathfrak{a}}\), there is an arch connecting the \(\left(2n-i\right)\)th and \(\left(2n-j\right)\)th points of \(e_{\mathfrak{a}\otimes\mathfrak{a}}\).  Thus if there is sequence of edges in the diagram \(\langle\xi_{\mathfrak{a}},\xi_{\mathfrak{b}}\rangle\) which connects the point at infinity of \(\xi_{\mathfrak{b}}\) back to itself, then, since it may not contain any of the through-strings of \(\xi_{\mathfrak{a}}\), there must be a string in the diagram \(e_{\mathfrak{a}\otimes\mathfrak{a}}\xi_{\mathfrak{b}}\) connecting its point at infinity back to itself.  Thus \(e_{\mathfrak{a}\otimes\mathfrak{a}}\xi_{\mathfrak{b}}=0\), and hence \(\langle\xi_{\mathfrak{a}},\xi_{\mathfrak{b}}\rangle=q^{-p}\langle\xi_{\mathfrak{a}},e_{\mathfrak{a}\otimes\mathfrak{a}}\xi_{\mathfrak{b}}\rangle^{\prime}=0\).

Likewise, if there is a string in the diagram \(\langle\xi_{\mathfrak{a}},\xi_{\mathfrak{b}}\rangle\) connecting the point at infinity of \(\xi_{\mathfrak{a}}\) back to itself, there must also be such a string in \(e_{\mathfrak{a}\otimes\mathfrak{a}}\xi_{\mathfrak{b}}\), so \(e_{\mathfrak{a}\otimes\mathfrak{a}}\xi_{\mathfrak{b}}=0\), and hence as above \(\langle\xi_{\mathfrak{a}},\xi_{\mathfrak{b}}\rangle^{\prime}=0\).

On the other hand, if the diagram \(\langle\xi_{\mathfrak{a}},\xi_{\mathfrak{b}}\rangle\) does not have any string connecting either infinity back to itself, then each must have the same number of through-strings.

Each of \(\mathfrak{a}\) and \(\mathfrak{b}\) can be expressed as a series of box additions on the minimal element \(\left(1,\ldots,1\right)\).  As box addition at \(i\) corresponds to multiplication by \(e_{i}\) in the first basis, \(\xi_{\mathfrak{a}}=e_{i_{1}}\cdots e_{i_{r}}\xi_{\left(1,\ldots,1\right)}\) and \(\xi_{\mathfrak{b}}=e_{j_{1}}\cdots e_{j_{s}}\xi_{\left(1\ldots,1\right)}\) for some sequence of indices \(i_{1},\ldots,i_{r}\) and \(j_{1},\ldots,j_{s}\).  As each of these terms is self-adjoint, \(\left(e_{i_{1}}\cdots e_{i_{r}}\right)^{*}=e_{i_{r}}\cdots e_{i_{1}}\), and hence
\begin{eqnarray*}\langle\xi_{\mathfrak{a}},\xi_{\mathfrak{b}}\rangle^{\prime}&=&\langle e_{i_{1}}\cdots e_{i_{r}}\xi_{\left(1,\ldots,1\right)},e_{j_{1}}\cdots e_{j_{s}}\xi_{\left(1,\ldots,1\right)}\rangle^{\prime}\\&=&\langle\xi_{\left(1,\ldots,1\right)},e_{i_{r}}\cdots e_{i_{1}}e_{j_{1}}\cdots e_{j_{s}}\xi_{\left(1,\ldots,1\right)}\rangle\end{eqnarray*}

As established, \(\xi_{\left(1,\ldots,1\right)}=q^{-p}e_{\left(1,\ldots,1\right)\otimes\left(1,\ldots,1\right)}\xi_{\left(1,\ldots,1\right)}\) and \(e_{\left(1\ldots,1\right)\otimes\left(1,\ldots,1\right)}\) is self-adjoint, so
\begin{eqnarray*}\langle\xi_{\mathfrak{a}},\xi_{\mathfrak{b}}\rangle^{\prime}&=&q^{-2p}\langle e_{\left(1,\ldots,1\right)\otimes\left(1,\ldots,1\right)}\xi_{\left(1,\ldots,1\right)},\\&&e_{j_{r}}\cdots e_{j_{1}}e_{i_{1}}\cdots e_{i_{s}}e_{\left(1,\ldots,1\right)\otimes\left(1,\ldots,1\right)}\xi_{\left(1,\ldots,1\right)}\rangle^{\prime}\\&=&q^{-2p}\langle\xi_{\left(1,\ldots,1\right)},\\&&e_{\left(1,\ldots,1\right)\otimes\left(1,\ldots,1\right)}e_{j_{r}}\cdots e_{j_{1}}e_{i_{1}}\cdots e_{i_{s}}e_{\left(1,\ldots,1\right)\otimes\left(1,\ldots,1\right)}\xi_{\left(1,\ldots,1\right)}\rangle^{\prime}\end{eqnarray*}
We will let \(e=e_{\left(1,\ldots,1\right)\otimes\left(1,\ldots,1\right)}e_{j_{r}}\cdots e_{j_{1}}e_{i_{1}}\cdots e_{i_{s}}e_{\left(1,\ldots,1\right)\otimes\left(1,\ldots,1\right)}\).

Each arch in \(e_{\left(1,\ldots,1\right)\otimes\left(1,\ldots,1\right)}\) connecting two points in \(\left\{1,\ldots,n\right\}\) must also appear in the product \(e\), and likewise any arch connecting two points in \(\left\{n+1,\ldots,2n\right\}\) must also appear in \(e\).  As the diagram \(\langle\xi_{\mathfrak{a}},\xi_{\mathfrak{b}}\rangle\) has no strings connecting either infinity back to itself, it is nonzero.  Since \(\langle\xi_{\mathfrak{a}},\xi_{\mathfrak{b}}\rangle=\langle\xi_{\left(1,\ldots,1\right)},e\xi_{\left(1,\ldots,1\right)}\rangle\) as well, \(e\neq 0\), and \(e\) and \(e^{*}\) have nonzero action on \(\xi_{\left(1,\ldots,1\right)}\).  Thus the points of \(e\) which are identified with through-strings of \(\left(1,\ldots,1\right)\) when this action is calculated may not be connected to points on their own side.  So \(e\) must be \(e_{\left(1,\ldots,1\right)\otimes\left(1,\ldots,1\right)}\) along with some closed loops.  So \(e=q^{m}e_{\left(1,\ldots,1\right)\otimes\left(1,\ldots,1\right)}\) for some \(m\).  Thus
\begin{eqnarray*}\langle\xi_{\mathfrak{a}},\xi_{\mathfrak{b}}\rangle^{\prime}&=&q^{-2p}\langle\xi_{\left(1,\ldots,1\right)},q^{m}e_{\left(1,\ldots,1\right)\otimes\left(1,\ldots,1\right)}\xi_{\left(1,\ldots,1\right)}\rangle\\&=&q^{m-p}\langle\xi_{\left(1,\ldots,1\right)},\xi_{\left(1,\ldots,1\right)}\rangle^{\prime}\end{eqnarray*}

Thus the value of the inner product on any two basis elements is fully determined by its value on \(\xi_{\left(1,\ldots,1\right)}\).  Extending this sesquilinearly to all of \(U\left(n\right)\), we can see that our original inner product is unique up to scalar multiples.
\end{proof}
\end{thm}

\subsection{Orthogonality Theorems}

We now show that the second basis, consisting of the \(\xi_{\mathfrak{a}}^{\prime}\), is orthogonal.

\begin{lem}
\label{slopelem}
If \(\mathfrak{a}\) has a slope at \(i\), then \(e_{i}\xi^{\prime}_{\mathfrak{a}}=0\).
\begin{proof}
We prove this by induction on the number of boxes in \(a\).

\medskip
\noindent
\textit{Base Case:\ }
In our base case, \(\mathfrak{a}=\left(1,\ldots,1\right)\).  A slope at position \(i\) corresponds to two opening brackets or two closing brackets in positions \(i\) and \(i+1\).  In the first case \(h_{i+1}\left(\mathfrak{a}\right)>1\), and in the second, \(h_{i-1}\left(\mathfrak{a}\right)>1\).  Since the base height inside all of the arches of \(\left(1,\ldots,1\right)\) is \(1\), the point of height greater than \(1\) cannot occur inside an arch.  As closing brackets always correspond to a point connected by an arch, the second possibility, that there are two closing brackets at positions \(i\) and \(i+1\), is excluded.  If there are two opening brackets at positions \(i\) and \(i+1\), then the one at position \(i+1\), which has depth greater than \(1\), must correspond to a through-string.  As a through-string cannot be contained within an arch, the opening bracket at position \(i\) must also correspond to a through-string.  So a slope at position \(i\) in \(\mathfrak{a}=\left(1,\ldots,1\right)\) must correspond to through-strings at positions \(i\) and \(i+1\).

As \(e_{1}\) has an arch which connects the \(i\)th and \(\left(i+1\right)\)th points, the two through-strings of \(\xi^{\prime}_{\mathfrak{a}}=\xi_{\left(1,\ldots,1\right)}\) at positions \(i\) and \(i+1\) are connected in the diagram \(e_{i}\xi_{1,\ldots,1}\).  So \(e_{i}\xi^{\prime}_{\mathfrak{a}}=e_{i}\xi_{\left(1,\ldots,1\right)}=0\).

\medskip
\noindent
\textit{Inductive Step:\ }
We now assume that \(\mathfrak{a}\) has at least one box, so \(\mathfrak{a}=\Diamond_{j}\left(\mathfrak{b}\right)\) for some \(j\) and some \(\mathfrak{b}\) which has fewer boxes than \(\mathfrak{a}\).  A slope may be either increasing or decreasing.  We deal with the case where the slope at \(i\) is an increasing slope first.

\smallskip
\textit{Case I: the slope is increasing:\ }
Since \(h_{i}\left(\mathfrak{a}\right)=h_{i-1}\left(\mathfrak{a}\right)+1\), \(\mathfrak{a}\) does not have a maximum at \(i-1\), and since \(h_{i+1}\left(\mathfrak{a}\right)=h_{i}\left(\mathfrak{a}\right)+1\), it does not have a maximum at \(i\).  So \(j\) is not \(i-1\) or \(i\).

If \(j\) is also not \(i+1\), then the operator associated with the box addition commutes with multiplication by \(e_{i}\).  So \(e_{i}\xi^{\prime}_{\mathfrak{a}}=e_{i}\left(e_{j}-\mu_{h_{j}\left(\mathfrak{a}\right)-1}\right)\xi^{\prime}_{\mathfrak{b}}=\left(e_{j}-\mu_{h_{j}\left(\mathfrak{a}\right)-1}\right)e_{i}\xi^{\prime}_{\mathfrak{b}}\).  \(\mathfrak{b}\) has fewer boxes than \(\mathfrak{a}\), and since the box addition at \(j\) only affects the height at \(j\), the height of \(\mathfrak{b}\) at \(i-1\), \(i\) and \(i+1\) is the same as that of \(\mathfrak{a}\), and hence \(\mathfrak{b}\) also has a slope at \(i\).  So, by our inductive hypothesis, \(e_{i}\xi^{\prime}_{\mathfrak{b}}=0\).  Thus \(e_{i}\xi^{\prime}_{\mathfrak{a}}=0\).

If, on the other hand, \(j=i+1\), then \(\xi^{\prime}_{\mathfrak{a}}=\left(e_{i+1}-\mu_{h_{i+1}\left(\mathfrak{a}\right)-1}\right)\xi^{\prime}_{\mathfrak{b}}\), where \(\mathfrak{a}=\Diamond_{i+1}\left(\mathfrak{b}\right)\).

\(\mathfrak{a}\) has an increasing slope at \(i\), so its \(i\)th step is an upward step.  As box addition at \(i+1\) only affects the \(\left(i+1\right)\)th and \(\left(i+2\right)\)th steps, the \(i\)th step of \(\mathfrak{b}\) is also an upward step.  Since \(\mathfrak{b}\) admits a box addition at \(i+1\), it has a minimum at \(i+1\), and hence its \(\left(i+1\right)\)th step is a downward step.  So \(\mathfrak{b}\) has a maximum at \(i\).  So \(\xi^{\prime}_{\mathfrak{b}}=\left(e_{i}-\mu_{h_{i}\left(\mathfrak{b}\right)-1}\right)\xi\) for some \(\xi\in U\left(n;p\right)\).

We show that this is true regardless of whether the maximum at \(i\) is the result of box addition or not.  If the maximum at \(i\) in \(\mathfrak{b}\) is not the result of a box addition, then it must be present in the minimal Dyck path \(\left(1,\ldots,1\right)\), and hence in any path \(\mathfrak{c}\) with \(\left(1,\ldots,1\right)\preceq\mathfrak{c}\preceq\mathfrak{b}\).  So \(\xi^{\prime}_{\mathfrak{b}}\) is a linear combination of \(\xi_{\mathfrak{c}}\), each of which has a maximum at \(i\).  This maximum consists of an upward step at \(i\) and a downward step at \(i+1\), that is, an opening bracket at \(i\) and a closing bracket at \(i+1\).  These brackets must be paired, and hence correspond to an arch connecting the \(i\)th point to the \(\left(i+1\right)\)th point.  \(e_{i}\) acts on such a \(\xi_{\mathfrak{c}}\) by multiplying it by \(q\), so \(e_{i}\xi^{\prime}_{\mathfrak{b}}=q\xi^{\prime}_{\mathfrak{b}}\), and hence \(\left(e_{i}-\mu_{0}\right)\left(\frac{1}{q}\xi^{\prime}_{\mathfrak{b}}\right)=\xi^{\prime}_{\mathfrak{b}}\).

Thus, \(e_{i}\xi^{\prime}_{\mathfrak{a}}=e_{i}\left(e_{i+1}-\mu_{h_{i+1}\left(\mathfrak{a}\right)-1}\right)\left(e_{i}-\mu_{h_{i}\left(\mathfrak{b}\right)-1}\right)\xi\) for some \(\xi\in U\left(n;k\right)\).  \(\mathfrak{a}\) has an increasing slope at \(i\), so \(h_{i+1}\left(\mathfrak{a}\right)=h_{i}\left(\mathfrak{a}\right)+1\).  As \(\mathfrak{a}=\Diamond_{i+1}\left(\mathfrak{b}\right)\), the height at \(i\) is not altered, so \(h_{i}\left(\mathfrak{b}\right)=h_{i}\left(\mathfrak{a}\right)\).  Expanding, we get that \(\left(e_{i}e_{i+1}e_{i}-\mu_{h_{i}\left(\mathfrak{a}\right)-1}e_{i}e_{i+1}-\mu_{h_{i}\left(\mathfrak{a}\right)}e_{i}^{2}+\mu_{h_{i}\left(\mathfrak{a}\right)}\mu_{h_{i}\left(\mathfrak{a}\right)-1}e_{i}\right)\xi\).  By the relations for the Temperley-Lieb algebra, \(e_{i}e_{i+1}e_{i}=e_{i}\) and \(e_{i}^{2}=qe_{i}\).  We can rearrange our terms to get \(\left(1-q\mu_{h_{i}\left(\mathfrak{a}\right)}+\mu_{h_{i}\left(\mathfrak{a}\right)}\mu_{h_{i}\left(\mathfrak{a}\right)-1}\right)e_{i}\xi-\mu_{h_{i}\left(\mathfrak{a}\right)-1}e_{i}e_{i+1}\xi\).  By the recurrence relation of the Chebyshev polynomials, \(\Delta_{h_{i}\left(\mathfrak{a}\right)}-q\Delta_{h_{i}\left(\mathfrak{a}\right)-1}+\Delta_{h_{i}\left(\mathfrak{a}\right)-2}=0\), so dividing through by \(\Delta_{h_{i}\left(\mathfrak{a}\right)}\), \(1-q\frac{\Delta_{h_{i}\left(\mathfrak{a}\right)-1}}{\Delta_{h_{i}\left(\mathfrak{a}\right)}}+\frac{\Delta_{h_{i}\left(\mathfrak{a}\right)-2}}{\Delta_{h_{i}\left(\mathfrak{a}\right)-1}}\frac{\Delta_{h_{i}\left(\mathfrak{a}\right)-1}}{\Delta_{h_{i}\left(\mathfrak{a}\right)}}=0\), so the first term vanishes.  So \(e_{i}\xi^{\prime}_{\mathfrak{a}}=-\mu_{h_{i}\left(\mathfrak{a}\right)-1}e_{i}e_{i+1}\xi\).

If the previously discussed maximum in \(\mathfrak{b}\) at \(i\) is not the result of box addition, then it is a maxiumum in the minimal path \(\left(1,\ldots,1\right)\), corresponding to a pairing of the \(i\)th and \(\left(i+1\right)\)th points.  As all arches in \(\left(1,\ldots,1\right)\) have depth \(1\), \(h_{i-1}\left(\mathfrak{b}\right)=0\), so \(h_{i}\left(\mathfrak{a}\right)=1\).  Thus \(\mu_{h_{i}\left(\mathfrak{a}\right)-1}=\frac{\Delta_{-1}}{\Delta_{0}}=0\), so the second term also vanishes.

If the maximum in \(\mathfrak{b}\) at \(i\) is the result of box addition, then \(\xi=\xi^{\prime}_{\mathfrak{c}}\) where \(\mathfrak{b}=\Diamond_{i}\left(\mathfrak{c}\right)\).  Since \(\mathfrak{c}\) admits box addition at \(i\), it must have a minimum at \(i\), and hence its \(\left(i+1\right)\)th step is an upward step.  Its \(\left(i+2\right)\)th step will be the same as that of \(\mathfrak{b}\), and since \(\mathfrak{b}\) admits a box addition at \(i+1\), this must also be an upward step.  So \(\mathfrak{c}\) has an increasing slope at \(i+1\).  Since \(\mathfrak{c}\) has two fewer boxes than \(\mathfrak{a}\), by our induction hypothesis, \(e_{i+1}\xi^{\prime}_{\mathfrak{c}}=0\).  So the second term in our expression for \(e_{i}\xi^{\prime}_{\mathfrak{a}}\) also vanishes.

\smallskip
\textit{Case II: the slope is decreasing:\ }
Likewise, if \(\mathfrak{a}\) has a decreasing slope at \(i\), then it does not have a maximum at \(i\) or \(i+1\).  If it also does not have a maximum at \(i-1\), then the last box addition must have been elsewhere, and hence commutes with \(e_{i}\).  So, as above, \(e_{i}\xi^{\prime}_{\mathfrak{a}}=0\).  If the last box addition was at \(i-1\), then \(\xi^{\prime}_{\mathfrak{a}}=\left(e_{i-1}-\mu_{h_{i-1}\left(\mathfrak{a}\right)-1}\right)\left(e_{i}-\mu_{h_{i}\left(\mathfrak{a}\right)-1}\right)\xi\) for some \(\xi\in U\left(n;k\right)\).  Then as before \(e_{i}\xi^{\prime}_{\mathfrak{a}}=\left(e_{i}e_{i-1}e_{i}-\mu_{h_{i}\left(\mathfrak{a}\right)-1}e_{i}e_{i-1}-\mu_{h_{i}\left(\mathfrak{a}\right)}e_{i}^{2}+\mu_{h_{i}\left(\mathfrak{a}\right)}\mu_{h_{i}\left(\mathfrak{a}\right)-1}e_{i}\right)\xi=\left(1-q\mu_{h_{i}\left(\mathfrak{a}\right)}+\mu_{h_{i}\left(\mathfrak{a}\right)}\mu_{h_{i}\left(\mathfrak{a}\right)-1}\right)e_{i}\xi-\mu_{h_{i}\left(\mathfrak{a}\right)-1}e_{i}e_{i-1}\xi\).  As above, by the recurrence relations of the Chebyshev polynomials, the first term vanishes.  Again, either \(h_{i}\left(\mathfrak{a}\right)=1\) so \(\mu_{h_{i}\left(\mathfrak{a}\right)-1}=0\), or \(\xi=\xi^{\prime}_{\mathfrak{c}}\) for some \(\mathfrak{c}\) with a slope at \(i-1\) so that \(e_{i-1}\xi^{\prime}_{\mathfrak{c}}=0\), so the second term vanishes.
\end{proof}
\end{lem}

Any two basis elements corresponding to diagrams with different numbers of through-strings are orthogonal, since every edge-path beginning on a through-string must eventually end on a through-string, and since one diagram has more through-strings, not all of the edge-paths beginning on one of those may end on a through-string of the other diagram.  Thus the subspaces \(U\left(n;p\right)\) are orthogonal.  Since each element \(\xi^{\prime}_{\mathfrak{a}}\) is a linear combination of \(\xi_{\mathfrak{b}}\in U\left(n;p\right)\) for some fixed \(p\), \(\xi^{\prime}_{\mathfrak{a}}\in U\left(n;p\right)\).  So we know that the new basis elements corresponding to diagrams with different numbers of through-strings are orthogonal.  We now show that \(\xi^{\prime}_{\mathfrak{a}}\) and \(\xi^{\prime}_{\mathfrak{b}}\) are orthogonal even if \(\mathfrak{a}\) and \(\mathfrak{b}\) have the same number of through-strings.

\begin{thm}
\label{orthothm}
If \(\mathfrak{a}\neq\mathfrak{b}\), then \(\langle\xi^{\prime}_{\mathfrak{a}},\xi^{\prime}_{\mathfrak{b}}\rangle=0\).
\begin{proof}
We prove this by induction on the number of boxes in \(\mathfrak{a}\).

\medskip
\noindent
\textit{Base Case}
In our base case, \(\mathfrak{a}=\left(1,\ldots,1\right)\).  We prove by induction on the number of boxes in \(\mathfrak{b}\).  Since \(\mathfrak{a}\neq\mathfrak{b}\), we can assume that \(\mathfrak{b}\) is the result of at least one box addition.  We consider first the case where \(\mathfrak{b}=\Diamond_{i}\left(\mathfrak{c}\right)\) where \(\mathfrak{a}\) does not have a minimum at \(i\), then the case where there is no such \(i\) and \(\mathfrak{c}\).

\smallskip
\textit{Case I: there is an \(i\) such that \(\mathfrak{b}=\Diamond_{i}\left(\mathfrak{c}\right)\) and \(\mathfrak{a}\) does not have a minimum at \(i\):\ }
In the first case, \(\langle\xi^{\prime}_{\mathfrak{a}},\xi^{\prime}_{\mathfrak{b}}\rangle=\langle\xi^{\prime}_{\mathfrak{a}},\left(e_{i}-\mu_{h_{i}\left(\mathfrak{b}\right)-1}\right)\xi^{\prime}_{\mathfrak{c}}\rangle=\langle e_{i}\xi^{\prime}_{\mathfrak{a}},\xi^{\prime}_{\mathfrak{c}}\rangle-\mu_{h_{i}\left(\mathfrak{b}\right)-1}\langle\xi^{\prime}_{\mathfrak{a}},\xi^{\prime}_{\mathfrak{b}}\rangle\).  As \(\mathfrak{a}\) does not have a minimum at \(i\), it has either a slope or a maximum.  By the previous lemma, \(e_{i}\xi^{\prime}_{\mathfrak{a}}=0\).  If it is a maximum, then as established in the same lemma, \(e_{i}\xi^{\prime}_{\mathfrak{a}}=q\xi^{\prime}_{\mathfrak{a}}\).  Either way, it is a scalar multiple of \(\xi^{\prime}_{\mathfrak{a}}\), so the first term, and hence the entire expression, is a scalar multiple of \(\langle\xi^{\prime}_{\mathfrak{a}},\xi^{\prime}_{\mathfrak{c}}\rangle\).

Since \(\mathfrak{a}\) does not have a minimum at \(i\) and hence does not admit a box addition at \(i\), and \(\mathfrak{c}\), \(\mathfrak{a}\neq\mathfrak{c}\).  As \(\mathfrak{c}\) has fewer boxes than \(\mathfrak{b}\), by our inductive hypothesis, we have that \(\langle\xi^{\prime}_{\mathfrak{a}},\xi^{\prime}_{\mathfrak{c}}\rangle=0\).  So the entire expression vanishes.

\smallskip
\textit{Case II: there is no such \(i\):\ }
In the second case, let \(i\) be the smallest index such that \(\mathfrak{b}=\Diamond_{i}\left(\mathfrak{c}\right)\) for some \(\mathfrak{c}\).  Then the \(i\)th step of \(\mathfrak{b}\) is an upward step and \(h_{i}\left(\mathfrak{b}\right)>1\).  If the \(\left(i-1\right)\)th step were a downward step, then \(h_{i-2}\left(\mathfrak{b}\right)>1\), so there must be a maximum at some index less than \(i\) with height greater than \(1\) which must hence be the result of a box addition.  Thus the \(\left(i-2\right)\)th step of \(\mathfrak{b}\) must also be an upward step, and hence \(\mathfrak{b}\) has a slope at \(i-1\).

Since \(\mathfrak{a}\) has a minimum at \(i\), its \(i\)th step is a downward step, which corresponds to a closing bracket, which is always part of an arch.  Since no arches in \(\mathfrak{a}=\left(1,\ldots,1\right)\) have depth greater than \(1\) and a through-string cannot be contained in an arch, there can be nothing inside this arch, so its opening bracket must be at \(i-1\).  So \(\mathfrak{a}\) must have a maximum at \(i-1\).  So, as established above, \(\xi^{\prime}_{\mathfrak{a}}=e_{i-1}\left(\frac{1}{q}\xi^{\prime}_{\mathfrak{a}}\right)\).

Thus, \(\langle\xi^{\prime}_{\mathfrak{a}},\xi^{\prime}_{\mathfrak{b}}\rangle=\langle\frac{1}{q}\xi^{\prime}_{\mathfrak{a}},e_{i-1}\xi^{\prime}_{\mathfrak{b}}\rangle\).  As \(\mathfrak{b}\) has a slope at \(i-1\), \(e_{i-1}\xi^{\prime}_{\mathfrak{b}}=0\).  So this entire quantity vanishes.

\medskip
\noindent
\textit{Inductive Step:\ }
In our inductive step, we can assume that \(\mathfrak{a}\) has at least one box, so \(\mathfrak{a}=\Diamond_{i}\left(\mathfrak{c}\right)\) for some \(i\) and some \(\mathfrak{c}\) with fewer boxes than \(\mathfrak{a}\), so our inductive hypothesis applies to \(\mathfrak{c}\).  We deal with the cases where \(\mathfrak{b}\) has a minimum, a slope, and a maximum at \(i\) separately.

\smallskip
\textit{Case I: \(\mathfrak{b}\) has a minimum at \(i\):\ }
Assume first that \(\mathfrak{b}\) has a minumum at \(i\).  Then \(\mathfrak{b}\) admits a box addition at \(i\), where \(\xi^{\prime}_{\Diamond_{i}\left(\mathfrak{a}\right)}=\left(e_{i}-\mu_{h_{i}\left(\mathfrak{b}\right)+1}\right)\xi^{\prime}_{\mathfrak{b}}\).  Then
\begin{eqnarray*}\langle\xi^{\prime}_{\mathfrak{a}},\xi^{\prime}_{\mathfrak{b}}\rangle&=&\langle\left(e_{i}-\mu_{h_{i}\left(\mathfrak{c}\right)+1}\right)\xi^{\prime}_{\mathfrak{c}},\xi^{\prime}_{\mathfrak{b}}\rangle\\&=&\langle\xi^{\prime}_{\mathfrak{c}},\left(e_{i}-\mu_{h_{i}\left(\mathfrak{b}\right)+1}\right)\xi^{\prime}_{\mathfrak{b}}\rangle-\langle\xi^{\prime}_{\mathfrak{c}},\left(\mu_{h_{i}\left(\mathfrak{c}\right)+1}-\mu_{h_{i}\left(\mathfrak{b}\right)+1}\right)\xi^{\prime}_{\mathfrak{b}}\rangle\\&=&\langle\xi^{\prime}_{\mathfrak{c}},\xi^{\prime}_{\Diamond_{i}\left(\mathfrak{b}\right)}\rangle-\left(\mu_{h_{i}\left(\mathfrak{c}\right)+1}-\mu_{h_{i}\left(\mathfrak{b}\right)+1}\right)\langle\xi^{\prime}_{\mathfrak{c}},\xi^{\prime}_{\mathfrak{b}}\rangle\end{eqnarray*}

As \(\mathfrak{c}\) admits box addition at \(i\), it has a minimum at \(i\), and since \(\Diamond_{i}\left(\mathfrak{b}\right)\) is the result of box addition at \(i\), it has a maximum at \(i\).  So \(\mathfrak{c}\neq\mathfrak{b}\), and hence \(\langle\xi^{\prime}_{\mathfrak{c}},\xi^{\prime}_{\Diamond_{i}\left(\mathfrak{b}\right)}\rangle=0\), so the first term vanishes.

If \(\mathfrak{c}\neq\mathfrak{b}\), then the second term vanishes as well.  If they are equal, then \(h_{i}\left(\mathfrak{c}\right)=h_{i}\left(\mathfrak{c}\right)\), so the two quotients of Chebyshev polynomials in the second term are equal, so their difference, and hence the second term, vanish.

\smallskip
\textit{Case II: \(\mathfrak{b}\) has a slope at \(i\):\ }
In the second case, \(\mathfrak{b}\) has a slope at \(i\).  Then
\begin{eqnarray*}\langle\xi^{\prime}_{\mathfrak{a}},\xi^{\prime}_{\mathfrak{b}}\rangle&=&\langle\left(e_{i}-\mu_{h_{i}\left(\mathfrak{a}\right)-1}\right)\xi^{\prime}_{\mathfrak{c}},\xi^{\prime}_{\mathfrak{b}}\rangle\\&=&\langle\xi^{\prime}_{\mathfrak{c}},e_{i}\xi^{\prime}_{\mathfrak{b}}\rangle-\mu_{h_{i}\left(\mathfrak{a}\right)-1}\langle\xi^{\prime}_{\mathfrak{c}},\xi^{\prime}_{\mathfrak{b}}\rangle\end{eqnarray*}
As \(\mathfrak{b}\) has as slope at \(i\), \(e_{i}\xi^{\prime}_{\mathfrak{b}}=0\), and hence the first term vanishes.  Since \(\mathfrak{c}\) has a minimum at \(i\) and \(\mathfrak{b}\) has a slope at \(i\), \(\mathfrak{c}\neq\mathfrak{b}\), so the second term also vanishes.

\smallskip
\textit{Case III: \(\mathfrak{b}\) has a maximum at \(i\):\ }
In the third case, \(\mathfrak{b}\) has a maximum at \(i\), so it can be expressed \(\left(e_{i}-\mu_{h_{i}\left(\mathfrak{b}\right)-1}\right)\xi\) for some \(\xi\in S\left(n;k\right)\).  As discussed above, this is true regardless of whether the maximum at \(i\) is the result of box addition or not.  If it is not, then \(\xi\) is a multiple of \(\xi^{\prime}_{\mathfrak{b}}\), and if it is, then \(\xi=\xi^{\prime}_{\mathfrak{d}}\) where \(\mathfrak{b}=\Diamond_{i}\left(\mathfrak{d}\right)\).  So \begin{eqnarray*}\langle\xi^{\prime}_{\mathfrak{a}},\xi^{\prime}_{\mathfrak{b}}\rangle&=&\langle\left(e_{i}-\mu_{h_{i}\left(\mathfrak{a}\right)-1}\right)\xi^{\prime}_{\mathfrak{c}},\left(e_{i}-\mu_{h_{i}\left(\mathfrak{b}\right)-1}\right)\xi\rangle\\&=&\langle\xi^{\prime}_{\mathfrak{c}},\left(qe_{i}-\mu_{h_{i}\left(\mathfrak{b}\right)-1}e_{i}-\mu_{h_{i}\left(\mathfrak{a}\right)-1}e_{i}+\mu_{h_{i}\left(\mathfrak{a}\right)-1}\mu_{h_{i}\left(\mathfrak{b}\right)-1}\right)\xi\rangle\\&=&\left(q-\mu_{h_{i}\left(\mathfrak{a}\right)-1}-\mu_{h_{i}\left(\mathfrak{b}\right)-1}\right)\langle\xi^{\prime}_{\mathfrak{c}},\left(e_{i}-\mu_{h_{i}\left(\mathfrak{b}\right)-1}\right)\xi\rangle\\&&+\left(\mu_{h_{i}\left(\mathfrak{a}\right)-1}\mu_{h_{i}\left(\mathfrak{b}\right)-1}+q\mu_{h_{i}\left(\mathfrak{b}\right)-1}-\mu_{h_{i}\left(\mathfrak{a}\right)-1}\mu_{h_{i}\left(\mathfrak{b}\right)-1}-\mu_{h_{i}\left(\mathfrak{b}\right)-1}^{2}\right)\\&&\langle\xi^{\prime}_{\mathfrak{c}},\xi\rangle\\&=&\left(q-\mu_{h_{i}\left(\mathfrak{a}\right)-1}-\mu_{h_{i}\left(\mathfrak{b}\right)-1}\right)\langle\xi^{\prime}_{\mathfrak{c}},\xi^{\prime}_{\mathfrak{b}}\rangle+\left(q\mu_{h_{i}\left(\mathfrak{b}\right)-1}-\mu_{h_{i}\left(\mathfrak{b}\right)-1}^{2}\right)\langle\xi^{\prime}_{\mathfrak{c}},\xi\rangle\end{eqnarray*}

As \(\mathfrak{c}\) has a minimum at \(i\) and \(\mathfrak{b}\) has a maximum at \(i\), \(\mathfrak{c}\neq\mathfrak{b}\), so the first term vanishes.  If the maximum at \(i\) in \(\mathfrak{b}\) is not the result of box addition, then \(\xi\) is a multiple of \(\xi^{\prime}_{\mathfrak{b}}\), so the second term also vanishes.  Otherwise, since \(\Diamond_{i}\left(\mathfrak{c}\right)\neq\Diamond_{i}\left(\mathfrak{d}\right)\), \(\mathfrak{c}\neq\mathfrak{d}\), so \(\langle\xi^{\prime}_{\mathfrak{c}},\xi^{\prime}_{\mathfrak{d}}\rangle=0\), and hence the second term also vanishes.
\end{proof}
\end{thm}

We get our first calculation of the square-norm of the second basis elements from Genauer and Stoltzfus \cite{2005math.....11003G} since Di Francesco \cite{MR1608551} and Cautis and Jackson \cite{MR1999738} normalize their basis elements.  However, we follow a method of proof more similar to the Di Francesco and Cautis and Jackson proofs.

\begin{thm}
\label{normthm}
If \(\mathfrak{a}=\left(a_{1},\ldots,a_{m}\right)\) then \(\langle\xi^{\prime}_{\mathfrak{a}},\xi^{\prime}_{\mathfrak{a}}\rangle=\frac{1}{\mu_{a_{1}}\cdots\mu_{a_{m}}}\).
\begin{proof}
We again prove this by induction on the number of boxes in \(\mathfrak{a}\).

\medskip
\noindent
\textit{Base Case:\ }
In our base case \(\xi^{\prime}_{\mathfrak{a}}=\xi_{\left(1,\ldots,1\right)}\).  In the diagram of the inner product \(\langle\xi_{\left(1,\ldots,1\right)},\xi_{\left(1,\ldots,1\right)}\rangle\), each of the \(p\) loops of one diagram \(\left(1,\ldots,1\right)\) is connected at both ends to the corresponding loop of the other diagram \(\left(1,\ldots,1\right)\), producing \(p\) closed loops and hence a factor of \(q^{p}\).  Likewise, each through-string of the upper \(\xi_{\left(1,\ldots,1\right)}\) is connected to the corresponding through-string of the lower \(\xi_{\left(1,\ldots,1\right)}\), producing a factor of \(1\).  So \(\langle\xi^{\prime}_{\mathfrak{a}},\xi^{\prime}_{\mathfrak{a}}\rangle=\langle\xi_{\left(1,\ldots,1\right)},\xi_{\left(1,\ldots,1\right)}\rangle=q^{p}=\frac{1}{\mu_{1}\cdots\mu_{1}}\).

\smallskip
\noindent
\textit{Inductive Step:\ }
For our induction step, we can assume that \(\mathfrak{a}=\left(a_{1},\ldots,a_{p}\right)\) has at least one box, so \(\mathfrak{a}=\Diamond_{i}\left(\mathfrak{b}\right)\) for some \(\mathfrak{b}=\left(b_{1},\ldots,b_{p}\right)\) with fewer boxes, and our induction hypothesis is that, \(\langle\xi^{\prime}_{\mathfrak{b}},\xi^{\prime}_{\mathfrak{b}}\rangle=\frac{1}{\mu_{b_{1}}\cdots\mu_{b_{p}}}\).

Thus,
\begin{eqnarray*}\langle\xi^{\prime}_{\mathfrak{a}},\xi^{\prime}_{\mathfrak{a}}\rangle&=&\langle\left(e_{i}-\mu_{h_{i}\left(\mathfrak{a}\right)-1}\right)\xi^{\prime}_{\mathfrak{b}},\left(e_{i}-\mu_{h_{i}\left(\mathfrak{a}\right)-1}\right)\xi^{\prime}_{\mathfrak{b}}\rangle\\&=&\langle\xi^{\prime}_{\mathfrak{b}},\left(qe_{i}-2\mu_{h_{i}\left(\mathfrak{a}\right)-1}e_{i}+\mu_{h_{i}\left(\mathfrak{a}\right)-1}^{2}\right)\xi^{\prime}_{\mathfrak{b}}\rangle\\&=&\left(q-2\mu_{h_{i}\left(\mathfrak{a}\right)-1}\right)\langle\xi^{\prime}_{\mathfrak{b}},\left(e_{i}-\mu_{h_{i}\left(\mathfrak{a}\right)-1}\right)\xi^{\prime}_{\mathfrak{b}}\rangle\\&&+\left(\mu_{h_{i}\left(\mathfrak{a}\right)-1}^{2}+q\mu_{h_{i}\left(\mathfrak{a}\right)-1}-2\mu_{h_{i}\left(\mathfrak{a}\right)-1}^{2}\right)\langle\xi^{\prime}_{\mathfrak{b}},\xi^{\prime}_{\mathfrak{b}}\rangle\\&=&\left(q-2\mu_{h_{i}\left(\mathfrak{a}\right)-1}\right)\langle\xi^{\prime}_{\mathfrak{b}},\xi^{\prime}_{\mathfrak{a}}\rangle+\mu_{h_{i}\left(\mathfrak{a}\right)-1}\left(q-\mu_{h_{i}\left(\mathfrak{a}\right)-1}\right)\langle\xi^{\prime}_{\mathfrak{b}},\xi^{\prime}_{\mathfrak{b}}\rangle\end{eqnarray*}
As \(\mathfrak{a}\neq\mathfrak{b}\), the first term vanishes, by the previous theorem.  By the recurrence relation for the Chebyshev polynomials, \(q-\mu_{h_{i}\left(\mathfrak{a}\right)-1}=\frac{q\Delta_{h_{i}\left(\mathfrak{a}\right)-1}-\Delta_{h_{i}\left(\mathfrak{a}\right)-2}}{\Delta_{h_{i}\left(\mathfrak{a}\right)-1}}=\frac{\Delta_{h_{i}\left(\mathfrak{a}\right)}}{\Delta_{h_{i}\left(\mathfrak{a}\right)-1}}=\frac{1}{\mu_{h_{i}\left(\mathfrak{a}\right)}}\).  Thus \(\langle\xi^{\prime}_{\mathfrak{a}},\xi^{\prime}_{\mathfrak{a}}\rangle=\frac{\mu_{h_{i}\left(\mathfrak{a}\right)-1}}{\mu_{h_{i}\left(\mathfrak{a}\right)}}\langle\xi^{\prime}_{\mathfrak{b}},\xi^{\prime}_{\mathfrak{b}}\rangle\).

We know that if \(\mathfrak{a}=\left(a_{1},\ldots,a_{p}\right)\) and \(\mathfrak{b}=\left(b_{1},\ldots,b_{p}\right)\), then \(\left(a_{1},\ldots,a_{p}\right)=\left(b_{1},\ldots,b_{j-1},b_{j}+1,b_{j+1},\ldots,b_{p}\right)\), where the decreasing step at \(i\) in \(\mathfrak{b}\) is its \(j\)th decreasing step, and the decreasing step at \(i+1\) in \(\mathfrak{a}\) is its \(j\)th decreasing step.  So \(b_{j}=h_{i-1}\left(\mathfrak{b}\right)\) and \(a_{j}=h_{i}\left(\mathfrak{a}\right)\).  Thus \(\langle\xi^{\prime}_{\mathfrak{a}},\xi^{\prime}_{\mathfrak{a}}\rangle=\frac{\mu_{b_{j}}}{\mu_{a_{j}}}\langle\xi^{\prime}_{\mathfrak{b}},\xi^{\prime}_{\mathfrak{b}}\rangle=\frac{\mu_{b_{j}}}{\mu_{a_{j}}}\frac{1}{\mu_{b_{1}}\cdots\mu_{b_{p}}}=\frac{1}{\mu_{a_{1}}\cdots\mu_{a_{p}}}\), as expected.
\end{proof}
\end{thm}

\begin{cor}
If \(\mathfrak{a}=\left(a_{1},\ldots,a_{p}\right)\) has maxima at \(i_{1},\ldots,i_{r}\) and minima at \(j_{1},\ldots,j_{s}\) (where \(\mathfrak{a}\) is thought to have a minimum at \(n\) if its last step is decreasing, but not to have a maximum at \(n\) if its last step is increasing), then \(\langle\xi^{\prime}_{\mathfrak{a}},\xi^{\prime}_{\mathfrak{a}}\rangle=\frac{\Delta_{h_{i_{1}}\left(\mathfrak{a}\right)}\cdots\Delta_{h_{i_{r}}\left(\mathfrak{a}\right)}}{\Delta_{h_{j_{1}}\left(\mathfrak{a}\right)}\cdots\Delta_{h_{i_{s}}\left(\mathfrak{a}\right)}}\).
\begin{proof}
We know that \(\langle\xi^{\prime}_{\left(1,\ldots,1\right)},\xi^{\prime}_{\left(1,\ldots,1\right)}\rangle=\langle\xi_{\left(1,\ldots,1\right)},\xi_{\left(1,\ldots,1\right)}\rangle=q^{p}\).  The height before each downward step is \(1\), so each of these must be a maximum (since otherwise the downward step must be preceded by at least one step, and if it is preceded by a downward step, then the height before that downward step is greater than the height before our original downward step, but we have no such heights in our restricted sequence), and hence we have \(p\) maxima with height \(1\), contributing a factor of \(q^{p}\) to the right-hand side of our equation.  Any minimum must have height \(0\), since it must follow one of these downward steps, each contributing a factor of \(\frac{1}{1}=1\) to the right-hand side of our equation.  So the right-hand side of the equation is \(q^{p}\).

Assume that the equation holds for \(\mathfrak{a}\).  As shown above, box addition at \(i\), where \(\mathfrak{a}\) has a minimum at \(i\), contributes a factor of \(\frac{\mu_{h_{i}\left(\mathfrak{a}\right)-1}}{\mu_{h_{i}\left(\mathfrak{a}\right)}}=\frac{\Delta_{h_{i}\left(\mathfrak{a}\right)}\Delta_{h_{i}\left(\mathfrak{a}\right)+2}}{\Delta_{h_{i}\left(\mathfrak{a}\right)+1}^{2}}\).  This turns the minimum \(\mathfrak{a}\) has at \(i\) into a maximum whose height is greater by two.  According to the right-hand side of the equation, we would expect to gain a factor of \(\Delta_{h_{i}\left(\mathfrak{a}\right)+2}\) for the new maximum and to lose a factor of \(\frac{1}{\Delta_{h_{i}\left(\mathfrak{a}\right)}}\) for the minimum we have lost, which accounts for the numerator of the expression.  Furthermore, changing the \(i\)th step to an upward step and the \(\left(i+1\right)\)th to a downward step can either create a minimum of height \(h_{i\pm 1}\left(\mathfrak{a}\right)=h_{i}\left(\mathfrak{a}\right)+1\) or turn a maximum of that height into a slope, on each side.  Either, according to the right-hand side of the equation, should give us a factor of \(\frac{1}{\Delta_{h_{i}\left(\mathfrak{a}\right)-1}}\).  This accounts for the denominator of the expression.  So the equation holds for \(\Diamond_{i}\left(\mathfrak{a}\right)\) as well.
\end{proof}
\end{cor}

\section{Explicit Isomorphism with a Direct Sum of Matrix Algebras}

We can now express the isomorphism between \(TL_{n}\) and a direct sum of matrix algebras exhibited in \cite{MR1999738} in terms of the representation of \(TL_{n}\) on \(U\left(n\right)\) and the bases we have just constructed.

We showed above that there is an obvious mapping between \(TL_{n}\) and \(U\left(2n;n\right)\), the algebra over noncrossing half-diagrams on \(2n\) points with no through-strings.  We can define a basis of \(e^{\prime}_{\mathfrak{p}}\in TL_{n}\), corresponding to the basis of \(\xi^{\prime}_{\mathfrak{a}}\in U\left(2n;n\right)\) under this mapping.  As the value of the inner product is preserved under this mapping, this is an orthogonal basis of \(TL_{n}\).  It can be shown further, as in \cite{MR1999738}, that \(e^{\prime}_{\mathfrak{a}\otimes\mathfrak{b}}e^{\prime}_{\mathfrak{c}\otimes\mathfrak{d}}=0\) if \(\mathfrak{b}\neq\mathfrak{c}\).  Here, we instead prove that \(e^{\prime}_{\mathfrak{a}\otimes\mathfrak{b}}\xi^{\prime}_{\mathfrak{c}}=0\) when \(\mathfrak{b}\neq\mathfrak{c}\).  (Note that with our abuse of notation, \(e^{\prime}_{\mathfrak{a}\otimes\mathfrak{b}}\) will not generally be equal to \(\xi^{\prime}_{\mathfrak{a}}\otimes\xi^{\prime}_{\mathfrak{b}}\).)

\begin{lem}
Let \(\mathfrak{a}\) and \(\mathfrak{b}\) be half-diagrams on \(n\) points with the same number of through-strings (so \(\mathfrak{a}\otimes\mathfrak{b}\) is defined) and let \(\mathfrak{c}\) be a half-diagram on \(n\) points where \(\mathfrak{b}\neq\mathfrak{c}\).  Then \(e^{\prime}_{\mathfrak{a}\otimes\mathfrak{b}}\xi^{\prime}_{\mathfrak{c}}=0\).
\begin{proof}
This proof is very similar to the previous proofs, where we transfer a box addition operator from one factor to the other.  We begin by showing that the action of \(e_{i}\in TL_{2n}\), for \(1,\leq i<n\), on \(\xi\in U\left(2n;n\right)\) is equivalent to left-multiplication of \(e_{i}\in TL_{n}\) on the \(e\in TL_{n}\) equivalent to \(\xi\), and likewise the action of \(e_{2n-i}\in TL_{2n}\), for \(1\leq i<n\), on \(\xi\) is equivalent to right-multiplication of \(e_{i}\in TL_{n}\) on this \(e\in TL_{n}\).  Then, since the box addition operators consist of an \(e_{i}\) or an \(e_{i}\) and a scalar multiple, box additions to \(\xi\in U\left(2n;n\right)\) not at \(n\) can be expressed in terms of left or right multiplication on the equivalent element of \(TL_{n}\).

\begin{figure}
\centering
\[\includegraphics{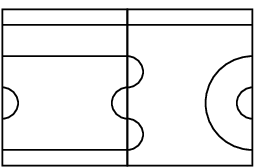}\simeq\includegraphics{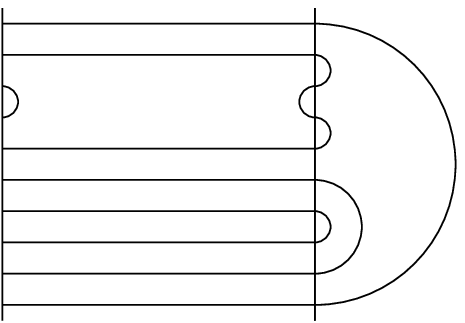}\]
\caption{Left multiplication by \(e_{3}\in TL_{5}\) is equivalent to the action of \(e_{3}\in TL_{10}\)}
\label{left}
\end{figure}

Let \(\mathfrak{p}\) be a full diagram on \(2n\) points, or alternately, a half-diagram on \(2n\) points with no through strings.  Then \(e_{i}e_{\mathfrak{p}}=e_{i}e_{\mathfrak{p}}\cdot 1\) and \(e_{\mathfrak{p}}e_{i}=1\cdot e_{\mathfrak{p}}e_{i}\).  We then wish to construct an element of \(TL_{2n}\) such that the same identifications of the same graphs occur, as shown in Figure \ref{left}.

In the first case, we construct a pairing on \(4n\) points by renumbering the total of \(4n\) points in \(e_{i}\) and \(1\).  We let the \(j\)th point of \(e_{i}\) be the \(j\)th point of the new diagram for \(1\leq j\leq n\), the \(\left(2n-j+1\right)\)th point of \(e_{i}\) be the \(\left(4n-j+1\right)\)th point of the new diagram for \(1\leq j\leq n\), the \(j\)th point of \(1\) be the \(\left(2n+j\right)\)th point of the new diagram for \(1\leq j\leq n\), and the \(\left(2n-j+1\right)\)th point of \(1\) be the \(\left(2n-j+1\right)\)th point of the new diagram, for \(1\leq j\leq n\).  Since each point is used exactly once, this is still a pairing.  We can see that the same points become points \(1\) through \(2n\) of the final diagram, and the same identifications of points are made.  Thus we get the same power of \(q\) and the same resulting diagram, but considered as an element of \(U\left(2n;n\right)\) rather than \(TL_{n}\).

Furthermore, we can see that since the \(i\)th point of \(e_{i}\in TL_{n}\) is connected to its \(\left(i+1\right)\)th, the \(i\)th point of the new diagram is connected to its \(\left(i+1\right)\)th, and since the \(\left(2n-i+1\right)\)th point of \(e_{i}\in TL_{n}\) is connected to its \(\left(2n-i\right)\)th, the \(\left(4n-i+1\right)\)th point of the new diagram is connected to its \(\left(4n-i\right)\)th.  If the \(j\)th point of the new diagram is not one of these points, it is connected to the \(\left(4n-j+1\right)\)th point.  We can recognize this element as \(e_{i}\in TL_{2n}\), so we know that it is a legitimate noncrossing diagram.

Thus, \(e_{i}e_{\mathfrak{p}}\) is equivalent to \(e_{i}\xi_{\mathfrak{p}}\).  We can extend this by linearity to all of \(TL_{n}\) and \(U\left(2n;n\right)\).

\begin{figure}
\centering
\[\includegraphics{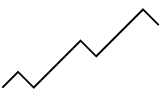}\otimes\includegraphics{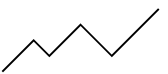}=\includegraphics{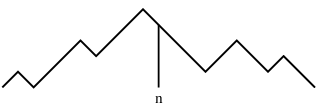}\]
\caption{\(\mathfrak{a}\), \(\mathfrak{b}\) and \(\mathfrak{a}\otimes\mathfrak{b}\) represented as Dyck paths: \(\mathfrak{a}\otimes\mathfrak{b}\) is \(\mathfrak{a}\) and the reflection of \(\mathfrak{b}\) in the line \(y=n\)}
\label{tensor}
\end{figure}

The points in \(\left\{1,\ldots,n\right\}\) in \(\mathfrak{a}\otimes\mathfrak{b}\) are paired as they are in \(\mathfrak{a}\), and those that were not paired in \(\mathfrak{a}\) are now paired with points in \(\left\{n+1,\ldots,2n\right\}\).  In either case, the \(i\)th bracket in \(\mathfrak{a}\otimes\mathfrak{b}\) is of the same type as the \(i\)th bracket in \(\mathfrak{a}\) for \(1\leq i\leq n\), and hence \(h_{i}\left(\mathfrak{a}\otimes\mathfrak{b}\right)=h_{i}\left(\mathfrak{a}\right)\), \(1\leq i\leq n\).  We show what this looks like when \(\mathfrak{a}\), \(\mathfrak{b}\) and \(\mathfrak{a}\otimes\mathfrak{b}\) are represented as Dyck paths in Figure \ref{tensor}.  So \(\Diamond_{i}\left(\mathfrak{a}\otimes\mathfrak{b}\right)=\Diamond_{i}\left(\mathfrak{a}\right)\otimes\mathfrak{b}\).  Thus \(\xi^{\prime}_{\Diamond_{i}\left(\mathfrak{a}\otimes\mathfrak{b}\right)}=\left(e_{i}-\mu_{h_{i}\left(\mathfrak{a}\otimes\mathfrak{b}\right)+1}\right)\xi^{\prime}_{\mathfrak{a}\otimes\mathfrak{b}}\), and hence \(e^{\prime}_{\Diamond_{i}\left(\mathfrak{a}\right)\otimes\mathfrak{b}}=\left(e_{i}-\mu_{h_{i}\left(\mathfrak{a}\right)+1}\right)e^{\prime}_{\mathfrak{a}\otimes\mathfrak{b}}\).  This is the appropriate operator for box addition on \(\mathfrak{a}\) at \(i\).

\begin{figure}
\centering
\[\includegraphics{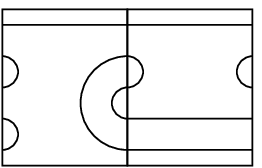}\simeq\includegraphics{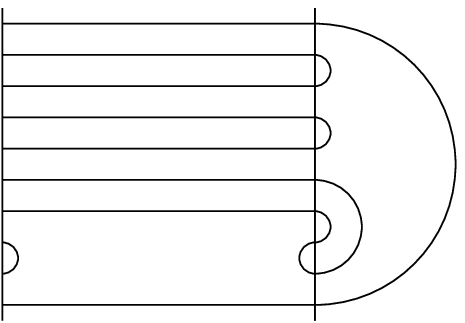}\]
\caption{Right multiplication by \(e_{2}\in TL_{5}\) is equivalent to the action of \(e_{8}\in TL_{10}\)}
\label{right}
\end{figure}

Likewise, \(e_{\mathfrak{p}}e_{i}=1\cdot e_{\mathfrak{p}}e_{i}\).  We can again construct an element of \(TL_{2n}\) which performs this multiplication in \(TL_{2n}\), shown in Figure \ref{right}.  We let the \(j\)th point of \(1\) be the \(j\)th point of the new diagram for \(1\leq j\leq n\), the \(\left(2n-j+1\right)\)th point of \(1\) be the \(\left(4n-j+1\right)\)th point of the new diagram, the \(j\)th point of \(e_{i}\) be the \(\left(2n+j\right)\)th point of the new diagram, and the \(\left(2n-j+1\right)\)th point of \(e_{i}\) be the \(\left(2n-j+1\right)\)th of the new diagram for \(1\leq j\leq n\).  Again, the same identifications are made, so we get the same power of \(q\) and the same resulting diagram.

The points \(i\) and \(i+1\) are paired in \(e_{i}\in TL_{n}\), so the points \(2n+i=4n-\left(2n-i\right)\) and \(2n+i+1=4n-\left(2n-i\right)+1\) are paired in the new diagram, and points \(2n-i+1\) and \(2n-i\) are paired in \(e_{i}\), so points \(2n-i+1\) and \(2n-i\) are paired in the new diagram.  For all other \(j\) in the new diagram, the \(j\)th point is connected to the \(\left(4n-j+1\right)\)th.  We can again recognize this pairing as \(e_{2n-i}\), a legitimate noncrossing pairing.

As above, we can extend this linearly, so \(e_{2n}\xi\) is equivalent to \(ee_{i}\) when \(\xi\in U\left(2n;n\right)\) is equivalent to \(e\in TL_{n}\).

\(\mathfrak{a}\otimes\mathfrak{b}\) is the full diagram whose \(i\)th point is the \(i\)th point of \(\mathfrak{a}\) for \(1\leq i\leq n\), and whose \(\left(2n-i+1\right)\)th point is the \(i\)th point of \(\mathfrak{b}\), for \(1\leq i\leq n\), with the through-strings of the two diagrams connected in the unique noncrossing way.  This reverses the order of the \(n\) points of \(\mathfrak{b}\), so the paired points corresponding to opening brackets now correspond to closing brackets, and those corresponding to closing brackets now correspond got opening brackets.  As the through-strings of \(\mathfrak{b}\) are now connected to points in \(\left\{1,\ldots,n\right\}\), they correspond to closing brackets.  So the \(\left(2n-i\right)\)th bracket of \(\mathfrak{a}\otimes\mathfrak{b}\) is the opposite type of bracket as the \(i\)th bracket of \(\mathfrak{b}\).  Thus if \(\mathfrak{b}\) has a maximum, minimum or slope at \(i\) (that is, the brackets at \(i\) and \(i+1\) are an opening bracket and a closing bracket, a closing bracket and an opening bracket, or two brackets of the same type), \(\mathfrak{a}\otimes\mathfrak{b}\) has a maximum, minimum or slope respectively at \(2n-i\) (since the brackets at \(2n-i\) and \(2n-i+1\) will be an opening bracket and a closing bracket, a closing bracket and an opening bracket, or the same type of bracket).

Since \(\mathfrak{a}\otimes\mathfrak{b}\) corresponds to a legitimate Dyck path, whose final point must be \(\left(2n,0\right)\), its height at \(2n-i\), for \(1\leq i\leq n\), must then be equal to the number of downward steps in the last \(i\) steps minus the number of upward steps in the last \(i\) steps, that is, the number by which the upward steps exceeds the downward steps among the first \(i\) of \(\mathfrak{b}\).  So \(h_{2n-i}\left(\mathfrak{a}\otimes\mathfrak{b}\right)=h_{i}\left(\mathfrak{b}\right)\).  This is shown in Figure \ref{tensor}.

Then we have \(\Diamond_{2n-i}\left(\mathfrak{a}\otimes\mathfrak{b}\right)=\mathfrak{a}\otimes\Diamond_{i}\left(\mathfrak{b}\right)\), and therefore \(\xi^{\prime}_{\Diamond_{2n-i}\left(\mathfrak{a}\otimes\mathfrak{b}\right)}=\left(e_{2n-i}-\right.\allowbreak\left.\mu_{h_{2n-i}\left(\mathfrak{a}\otimes\mathfrak{b}\right)+1}\right)\xi^{\prime}_{\mathfrak{a}\otimes\mathfrak{b}}\) is equivalent to \(e^{\prime}_{\mathfrak{a}\otimes\Diamond_{i}\left(\mathfrak{b}\right)}=e^{\prime}_{\mathfrak{a}\otimes\mathfrak{b}}\left(e_{i}-\mu_{h_{i}\left(\mathfrak{b}\right)+1}\right)\).

If \(\mathfrak{b}\) has a slope at \(i\), then \(\mathfrak{a}\otimes\mathfrak{b}\) has a slope at \(2n-i\), so \(e^{\prime}_{\mathfrak{a}\otimes\mathfrak{b}}e_{i}=e_{2n-i}\xi^{\prime}_{\mathfrak{a}\otimes\mathfrak{b}}=0\), so we have a version of Lemma \ref{slopelem} which we can apply in this situation.  We can then establish the base case of our induction on the number of boxes of \(\mathfrak{c}\) (which needs to be extended to the case where \(\mathfrak{b}\) and \(\mathfrak{c}\) are both minimal elements, but with different numbers of through-strings, which we do below).  The case where \(\mathfrak{b}\) has at least one box but \(\mathfrak{c}\) does not depends only on Lemma \ref{slopelem}, so the rest of the base case applies here as well.  The calculations in the three cases of the inductive step are identical as well.

\medskip
\noindent
\textit{Base Case:\ }
We first show that if \(\mathfrak{b}\) and \(\mathfrak{c}\) are both minimal elements in their respective orderings, but have different numbers of through-strings, then the statement holds.  Let \(p_{1}\) be the number of pairs in \(\mathfrak{b}\), and let \(p_{2}\) be the number of pairs in \(\mathfrak{c}\).  We consider the case where \(p_{1}<p_{2}\) and the case \(p_{1}>p_{2}\) separately.

\smallskip
\noindent
\textit{Case I: \(p_{1}<p_{2}\):\ }
\(\mathfrak{c}=\left(1,\ldots,1\right)\), so none of its \(p_{2}\) loops may be underneath another loop or through-string.  Thus each loop must close immediately, and the \(p_{2}\) loops must precede all of the \(n-2p_{2}\) through-strings.  So it has a loop from \(2i-1\) to \(2i\) for all \(i\), \(1\leq i\leq p_{2}\).  Thus \(e_{p_{2}}\xi^{\prime}_{\mathfrak{c}}=e_{p_{2}}\xi_{\mathfrak{c}}=q\xi_{\mathfrak{c}}\).  So \(e^{\prime}_{\mathfrak{a}\otimes\mathfrak{b}}\xi^{\prime}_{\mathfrak{c}}=\frac{1}{q}e^{\prime}_{\mathfrak{a}\otimes\mathfrak{b}}e_{p_{2}}\xi^{\prime}_{\mathfrak{c}}\).

\(\mathfrak{b}=\left(1,\ldots,1\right)\) as well, so again, none of its arches may occur underneath another arch or a through-string, so none of its \(n-2p_{1}\) through-strings must occur before the last \(n-2p_{1}\) points.  Thus \(\mathfrak{b}\) ends with \(n-2p_{1}\) upward steps, and hence has a slope at \(i\) for all \(i\), \(2p_{1}+1\leq i<n\), and hence \(\mathfrak{a}\otimes\mathfrak{b}\) has a slope at \(2n-i\) for all such \(i\).

Right multiplication of \(e^{\prime}_{\mathfrak{a}\otimes\mathfrak{b}}\) by \(e_{i}\) is equivalent to action of \(e_{2n-i}\) on the equivalent element of \(U\left(2n;n\right)\), \(\xi_{\mathfrak{a}\otimes\mathfrak{b}}\), and \(e_{2n-p_{2}}\xi^{\prime}_{\mathfrak{a}\otimes\mathfrak{b}}=0\).  So the entire expression vanishes, and \(e^{\prime}_{\mathfrak{a}\otimes\mathfrak{b}}\xi^{\prime}_{\mathfrak{c}}=0\).

\smallskip
\noindent
\textit{Case II: \(p_{1}>p_{2}\):\ }
\(e^{\prime}_{\mathfrak{a}\otimes\mathfrak{b}}\) is a linear combination of \(e_{\mathfrak{p}}\) with \(\mathfrak{p}\preceq\mathfrak{a}\otimes\mathfrak{b}\).  Thus \(h_{n}\left(\mathfrak{p}\right)\leq h_{n}\left(\mathfrak{a}\otimes\mathfrak{b}\right)\).  \(h_{n}\left(\mathfrak{p}\right)\) is the number by which the opening brackets exceed the closing brackets in the first \(n\) brackets of \(\mathfrak{p}\), that is, the number of points in \(\left\{1,\ldots,n\right\}\) which must be paired with a point in \(\left\{n+1,\ldots,2n\right\}\).  So any \(\mathfrak{p}\) with \(\mathfrak{p}\preceq\mathfrak{a}\otimes\mathfrak{b}\) may have at most as many through-strings as \(\mathfrak{a}\otimes\mathfrak{b}\).

Thus, \(e^{\prime}_{\mathfrak{a}\otimes\mathfrak{b}}\) is a linear combination of \(e_{\mathfrak{p}}\) where \(\mathfrak{p}\) has \(n-2p\) through-strings, for \(p>p_{2}\).  So \(\mathfrak{c}\) must have at least two more through-strings than \(\mathfrak{p}\), and hence when we concatenate \(e_{\mathfrak{p}}\) with \(\xi^{\prime}_{\mathfrak{c}}=\xi_{\mathfrak{c}}\) in the usual way, there must be at least two through-strings in \(\mathfrak{c}\) which are not connected in a sequence of edges to a through string of \(\mathfrak{p}\).  Thus they must be connected to each other, so \(e_{\mathfrak{p}}\xi^{\prime}_{\mathfrak{c}}=0\).  So \(e^{\prime}_{\mathfrak{a}\otimes\mathfrak{b}}\xi^{\prime}_{\mathfrak{p}}=0\).
\end{proof}
\end{lem}

We can also adapt our proof of the formula for \(\langle\mathfrak{a},\mathfrak{a}\rangle\) to the situation where \(\mathfrak{b}=\mathfrak{c}\).

\begin{lem}
Let \(\mathfrak{a}\) and \(\mathfrak{b}\) be half-diagrams on \(n\) points with the same number of through-strings.  Then \(e^{\prime}_{\mathfrak{a}\otimes\mathfrak{b}}\xi^{\prime}_{\mathfrak{b}}=\langle\xi^{\prime}_{\mathfrak{b}},\xi^{\prime}_{\mathfrak{b}}\rangle\xi^{\prime}_{\mathfrak{a}}\).
\begin{proof}
If we can express \(\mathfrak{b}\) as a box addition, we can perform the same calculations as in Theorem \ref{normthm} to get the same recurrence relations.  All we need to show is that, in the base case, where \(\mathfrak{b}\) is a minimal element, \(e^{\prime}_{\mathfrak{a}\otimes\mathfrak{b}}=\langle\xi^{\prime}_{\mathfrak{b}},\xi^{\prime}_{\mathfrak{b}}\rangle\xi^{\prime}_{\mathfrak{a}}\).  Then we can express \(\langle\xi^{\prime}_{\Diamond_{i}\left(\mathfrak{b}\right)},\xi^{\prime}_{\Diamond_{i}\left(\mathfrak{b}\right)}\rangle\) as a scalar multiple of \(\langle\xi^{\prime}_{\mathfrak{b}}\xi^{\prime}_{\mathfrak{b}}\rangle\), since all other cross-terms vanish.  We can likewise express \(e^{\prime}_{\mathfrak{a}\otimes\Diamond_{i}\left(\mathfrak{b}\right)}\xi^{\prime}_{\Diamond_{i}\left(\mathfrak{b}\right)}\) as a scalar mutliple of \(e^{\prime}_{\mathfrak{a}\otimes\mathfrak{b}}\xi^{\prime}_{\mathfrak{b}}\), since again all other cross-terms vanish.

\medskip
\noindent
\textit{Base Case:\ }
Let \(\mathfrak{b}\) be a minimal element, so \(\xi^{\prime}_{\mathfrak{b}}=\xi_{\mathfrak{b}}\).  \(\mathfrak{a}\) can be expressed as a series of box additions on \(\mathfrak{b}\) (since it has the same number of through-strings), so \(e^{\prime}_{\mathfrak{a}\otimes\mathfrak{b}}\) can be expressed as \(e^{\prime}_{\mathfrak{b}\otimes\mathfrak{b}}\) left-multiplied by the factors corresponding to these box additions.

\(e^{\prime}_{\mathfrak{b}\otimes\mathfrak{b}}\) is a linear combination of \(e_{\mathfrak{c}\otimes\mathfrak{d}}\) with \(\mathfrak{c}\otimes\mathfrak{d}\preceq\mathfrak{b}\otimes\mathfrak{b}\).  The coefficient of \(e_{\mathfrak{b}\otimes\mathfrak{b}}\) in this is \(1\), as shown above.  For any other \(\mathfrak{c}\otimes\mathfrak{d}\), \(h_{i}\left(\mathfrak{c}\otimes\mathfrak{d}\right)<h_{i}\left(\mathfrak{b}\otimes\mathfrak{b}\right)\) for some \(i\).  Since \(\mathfrak{b}\) is minimal, we cannot find a \(\mathfrak{c}\prec\mathfrak{b}\) or \(\mathfrak{d}\prec\mathfrak{b}\), so \(\mathfrak{c}\) and \(\mathfrak{d}\) must not have the same number of through-strings as \(\mathfrak{b}\).  Since the number of through-strings of \(\mathfrak{c}\) and \(\mathfrak{d}\) is \(h_{n}\left(\mathfrak{c}\otimes\mathfrak{d}\right)\), and \(h_{n}\left(\mathfrak{c}\otimes\mathfrak{d}\right)\leq h_{n}\left(\mathfrak{b}\otimes\mathfrak{b}\right)\), the \(\mathfrak{c}\otimes\mathfrak{d}\) must have fewer through strings.  Then, as shown above, \(e_{\mathfrak{c}\otimes\mathfrak{d}}\xi_{\mathfrak{b}}=0\).  So these terms do not contribute.

We also know from above that \(e_{\mathfrak{b}\otimes\mathfrak{b}}\xi_{\mathfrak{b}}=q^{p}\xi_{\mathfrak{b}}\), where \(p\) is the number of loops in \(\mathfrak{b}\).  So \(e^{\prime}_{\mathfrak{b}\otimes\mathfrak{b}}\xi^{\prime}_{\mathfrak{b}}=q^{p}\xi^{\prime}_{\mathfrak{b}}\).  Since \(e^{\prime}_{\mathfrak{a}\otimes\mathfrak{b}}\xi^{\prime}_{\mathfrak{b}}\) is \(e^{\prime}_{\mathfrak{b}\otimes\mathfrak{b}}\xi^{\prime}_{\mathfrak{b}}\) left-multiplied by the box addition operators which would give us \(\xi^{\prime}_{\mathfrak{a}}\) from \(\xi^{\prime}_{\mathfrak{b}}\), \(e^{\prime}_{\mathfrak{a}\otimes\mathfrak{b}}\xi^{\prime}_{\mathfrak{b}}=q^{p}\xi^{\prime}_{\mathfrak{b}}=\langle\xi^{\prime}_{\mathfrak{b}},\xi^{\prime}_{\mathfrak{b}}\rangle\xi^{\prime}_{\mathfrak{a}}\), as desired.
\end{proof}
\end{lem}

If we let the \(\xi^{\prime}_{\mathfrak{a}}\) be a basis for each \(U\left(n;p\right)\), then each \(e^{\prime}_{\mathfrak{p}}\) acts nontrivially on exactly one basis element in one of the \(U\left(n;p\right)\), which it maps to a scalar mutliple of another basis element.  So in this representation, each \(e^{\prime}_{\mathfrak{p}}\) is a nonzero multiple of an elementary matrix.  Thus we have an explicit isomorphism with the direct sum of matrix algebras \(M_{c_{n,0}\times c_{n,0}}\left(\mathbb{C}\right)\oplus\ldots\oplus M_{c_{n,2\lfloor\frac{n}{2}\rfloor}\times c_{n,2\lfloor\frac{n}{2}\rfloor}}\left(\mathbb{C}\right)\).

However, if we want this to be a \({}^{*}\)-homomorphism, we need to normalize the basis elements of the \(U\left(2n;n\right)\).  As \(\langle\xi^{\prime}_{\mathfrak{a}},\xi^{\prime}_{\mathfrak{a}}\rangle\) is not necessarily nonnegative, we cannot actually normalize the basis.  However, we can choose a square root of this quantity and take our basis to be \(\frac{\xi^{\prime}_{\mathfrak{a}}}{\sqrt{\langle\xi^{\prime}_{\mathfrak{a}},\xi^{\prime}_{\mathfrak{a}}\rangle}}\), with the convention that if the inner product is negative, we take the positive imaginary squareroot.

Then \(e^{\prime}_{\mathfrak{a}\otimes\mathfrak{b}}\) acts nontrivially only on normalized basis element \(\frac{\xi^{\prime}_{\mathfrak{b}}}{\sqrt{\langle\xi^{\prime}_{\mathfrak{b}},\xi^{\prime}_{\mathfrak{b}}\rangle}}\), and \(e^{\prime}_{\mathfrak{a}\otimes\mathfrak{b}}\frac{\xi^{\prime}_{\mathfrak{b}}}{\sqrt{\langle\xi^{\prime}_{\mathfrak{b}},\xi^{\prime}_{\mathfrak{b}}}}=\sqrt{\langle\xi^{\prime}_{\mathfrak{b}},\xi^{\prime}_{\mathfrak{b}}\rangle}\xi^{\prime}_{\mathfrak{a}}=\sqrt{\langle\xi^{\prime}_{\mathfrak{a}},\xi^{\prime}_{\mathfrak{a}}\rangle}\sqrt{\langle\xi^{\prime}_{\mathfrak{b}},\xi^{\prime}_{\mathfrak{b}}\rangle}\frac{\xi^{\prime}_{\mathfrak{a}}}{\sqrt{\langle\xi^{\prime}_{\mathfrak{a}},\xi^{\prime}_{\mathfrak{a}}\rangle}}\).  So it is a scalar multiple of the elementary matrix indexed by \(\mathfrak{a}\) and \(\mathfrak{b}\).

We now consider \(\left(e^{\prime}_{\mathfrak{a}\otimes\mathfrak{b}}\right)^{*}\).  This is constructed by relabelling the points in each of the component diagrams from \(i\) to \(2n-i+1\).  We can do this by relabelling the points of the minimal element from \(i\) to \(2n-i+1\), and by relabelling the points in any diagram in \(TL_{2n}\) which acts on it from \(i\) to \(2n-i+1\) for \(1\leq i\leq 2n\) and from \(2n+i\) to \(4n-i\) for \(1\leq i\leq 2n\).

No arch of the minimal element \(\left(1,\ldots,1\right)\) occurs underneath another arch.  When we relabel the points, no arch of the new diagram is underneath another arch, so it must also be \(\left(1,\ldots,1\right)\), so it is still the minimal element.  If the minimal element is expressed \(\mathfrak{a}\otimes\mathfrak{b}\), then \(\mathfrak{a}\) and \(\mathfrak{b}\) are both minimal half-diagrams with the same number of through-strings, so \(\mathfrak{a}=\mathfrak{b}\), and hence \(\mathfrak{a}\otimes\mathfrak{b}=\mathfrak{b}\otimes\mathfrak{a}\).  So in this case \(\left(e^{\prime}_{\mathfrak{a}\otimes\mathfrak{b}}\right)^{*}=e^{\prime}_{\mathfrak{b}\otimes\mathfrak{a}}\).

When we relabel the points of \(e_{i}\) for some \(i\), it becomes \(e_{2n-i}\).  So if \(\left(e^{\prime}_{\mathfrak{a}\otimes\mathfrak{b}}\right)^{*}=e^{\prime}_{\mathfrak{b}\otimes\mathfrak{a}}\), then inductively \(\left(e^{\prime}_{\Diamond_{i}\left(\mathfrak{a}\otimes\mathfrak{b}\right)}\right)^{*}=\left(\left(e_{i}-\mu_{h_{i}\left(\mathfrak{a}\otimes\mathfrak{b}\right)+1}\right)e^{\prime}_{\mathfrak{a}\otimes\mathfrak{b}}\right)^{*}=\left(e_{2n-i}-\mu_{h_{2n-i}\left(\mathfrak{b}\otimes\mathfrak{a}\right)+1}\right)e^{\prime}_{\mathfrak{b}\otimes\mathfrak{a}}=e^{\prime}_{\Diamond_{2n-i}\left(\mathfrak{b}\otimes\mathfrak{a}\right)}\).

Thus, by induction, \(\left(e^{\prime}_{\mathfrak{a}\otimes\mathfrak{b}}\right)^{*}=e^{\prime}_{\mathfrak{b}\otimes\mathfrak{a}}\).  This element acts nontrivially only on basis element \(\frac{\xi^{\prime}_{\mathfrak{a}}}{\sqrt{\langle\xi^{\prime}_{\mathfrak{a}},\xi^{\prime}_{\mathfrak{a}}\rangle}}\), which it maps to \(\sqrt{\langle\xi^{\prime}_{\mathfrak{a}},\xi^{\prime}_{\mathfrak{a}}\rangle}\sqrt{\langle\xi^{\prime}_{\mathfrak{b}},\xi^{\prime}_{\mathfrak{b}}\rangle}\frac{\xi^{\prime}_{\mathfrak{b}}}{\sqrt{\langle\xi^{\prime}_{\mathfrak{b}},\xi^{\prime}_{\mathfrak{b}}\rangle}}\).  So \(\left(e^{\prime}_{\mathfrak{a}\otimes\mathfrak{b}}\right)^{*}\) is the same scalar mutliple of the elementary matrix indexed by \(\mathfrak{b}\) and \(\mathfrak{a}\) as \(e^{\prime}_{\mathfrak{a}\otimes\mathfrak{b}}\) is of the elementary matrix indexed by \(\mathfrak{a}\) and \(\mathfrak{b}\).

Thus, in this representation, the \(e^{\prime}_{\mathfrak{p}}\) are scalar multiples of the elementary matrices, and this isomorphism is a \({}^{*}\)-homomorphism.

\bibliography{Temperley-Lieb}
\bibliographystyle{plain}

\end{document}